\documentclass[10pt,reqno]{amsart}

\usepackage[utf8]{inputenc}
\usepackage[T1]{fontenc}
\usepackage{lmodern}
\usepackage{microtype}
\usepackage{amsmath,amssymb,mathtools}
\usepackage{enumitem}
\usepackage{amsthm}
\usepackage{bm}
\usepackage{relsize}
\usepackage{todonotes}
\usepackage{comment}

\DeclareMathOperator{\arcsec}{arcsec}
\DeclareMathOperator{\arcosh}{arcosh}
\newcommand{\normsq}[1]{\left\lVert#1\right\rVert^2}
\newcommand{\norm}[1]{\left\lVert#1\right\rVert}
\newcommand{\abs}[1]{\lvert#1\rvert}
\newcommand{\abssq}[1]{\lvert#1\rvert^2}
\newcommand{\grad}{\nabla}
\newcommand{\laplacian}{\partial^2_{x}}

\newcommand{\defeq}{\mathrel{\mathop:}=}
\newcommand{\lintegral}[2]{\int_{#1}^{#2}}
\newcommand{\fintegral}[1]{\int_{#1}}
\makeatletter
\def\Xint#1{\mathchoice
  {\XXint\displaystyle\textstyle{#1}}%
  {\XXint\textstyle\scriptstyle{#1}}%
  {\XXint\scriptstyle\scriptscriptstyle{#1}}%
  {\XXint\scriptscriptstyle\scriptscriptstyle{#1}}%
  \!\int}
\def\XXint#1#2#3{{\setbox0=\hbox{$#1{#2#3}{\int}$}%
  \vcenter{\hbox{$#2#3$}}\kern-.5\wd0}}
\providecommand{\fint}{\Xint-}
\makeatother
\newcommand{\bint}[1]{\fint_{#1}}
\newcommand{\diff}[1]{\partial_{#1}}

\renewcommand{\phi}{\varphi}
\newcommand{\support}{\text{supp }}
\newcommand{\ip}[2]{\langle #1, #2 \rangle}

\newcommand{\bbN}{\mathbb{N}}
\providecommand{\linfty}[1]{\norm{#1}_{\infty}}

\providecommand{\btk}{T_k}

\providecommand{\wto}{\rightharpoonup}
\newcommand{\cembed}{\overset{c}{\hookrightarrow}}
\providecommand{\dist}{\text{dist}}

\DeclareMathOperator*{\esssup}{ess\,sup}

\usepackage{graphicx}
\usepackage{hyperref}
\hypersetup{colorlinks=true,linkcolor=blue,citecolor=blue,urlcolor=blue}
\usepackage{cleveref}

\newcommand{\bbR}{\mathbb{R}}

\newcommand{\bbrd}{\mathbb{R}^d}

\newcommand{\mc}[1]{\mathcal{#1}}

\numberwithin{equation}{section}
\theoremstyle{plain}
\newtheorem{theorem}{Theorem}
\numberwithin{theorem}{section}
\newtheorem*{theorem*}{Theorem}
\newtheorem{proposition}[theorem]{Proposition}
\newtheorem{lemma}[theorem]{Lemma}
\newtheorem{corollary}[theorem]{Corollary}
\theoremstyle{definition}
\newtheorem{definition}[theorem]{Definition}

\theoremstyle{remark}
\newtheorem{remark}[theorem]{Remark}

\title[Existence and Smoothing for a Degenerate Diffusion]{Existence and Smoothing for a Nondivergence-Form Degenerate Diffusion from Plasma-Wave Theory}

\author{William Porteous}
\address{Department of Mathematics, The University of Texas at Austin, 2515 Speedway, PMA 8.100, Austin, TX 78712, USA}
\email{wporteous@utexas.edu}

\author{Irene M. Gamba}
\address{Oden Institute for Computational Engineering and Sciences, The University of Texas at Austin, 201 E 24th Street, POB 5.102, Austin, TX 78712, USA}
\email{gamba@math.utexas.edu}

\author{Kun Huang}
\address{Department of Mathematics, Virginia Tech, 225 Stanger Street, 460 McBryde Hall, Blacksburg, VA 24061-1026, USA}
\email{kunhuang@vt.edu}


\begin{document}

\begin{abstract}
We prove existence and positive-time smoothing of weak solutions to
the degenerate parabolic Cauchy--Dirichlet problem
$\partial_t u = \rho_\lambda(x)u\,\partial_x^2 u
+ \rho_\lambda(x)g(x)u$ on the half-line, where $\rho_\lambda$
vanishes at the boundary and grows at infinity. 
This scalar problem arises from a formal reduction of the one-dimensional
quasilinear system for plasma waves, which couples an electron distribution to
a wave-energy spectrum: we develop a solution theory for the scalar problem that lifts to this system.
For the scalar problem, we construct nonnegative weak solutions under
weighted integrability assumptions on the initial data and boundedness of the reaction coefficient $\rho_\lambda g$, 
admitting unbounded and discontinuous data pertinent to the physical model.
We establish a one-sided bound on $\partial_t u$ which yields
positive-time smoothing: solutions become locally H\"older
continuous in both space and time and locally Lipschitz continuous in space. Explicit examples show that this spatial regularity is sharp. 
Under natural assumptions on the initial data, our scalar
solutions lift to weak solutions of the particle--wave system that
preserve nonnegativity, conserve particle mass, smooth initial
particle measures to locally bounded densities at positive times,
and exhibit quantitative decay of the wave spectrum at large wavenumber.
\end{abstract}

\maketitle

\section{Introduction}\label{section:introduction}
This work addresses the nondivergence-form degenerate parabolic
Cauchy--Dirichlet problem~\eqref{eq:classical_form} below. We first review its origins in 
plasma-wave theory, then present the analytical framework and main results. This problem arises from the
\emph{quasilinear theory} of plasma waves, introduced by Drummond and
Pines~\cite{drummond} and Vedenov et al.~\cite{vedenov} to model the
relaxation of kinetic instabilities in collisionless plasma. This theory expands the Vlasov--Poisson system about a spatially
homogeneous background electron distribution. At first order,
this background distribution is held fixed and the spatially
oscillating perturbation obeys Landau's linear
theory~\cite{landau1946}; quasilinear theory then supplies a closure at second order,
describing the evolution of the background due to the
perturbation~\cite[Ch.~23]{thorne_blandford}. In one dimension,
this gives the system~\eqref{eq:system} with dispersion
relation $\omega(k)$:
\begin{equation}\label{eq:system}
\begin{split}
 &\partial_t f(t,p) = \partial_p \left(
  \left[ \int_{\bbR} W(t,k)\, \delta_{0}\bigl( \omega(k) - pk
  \bigr) dk \right] \partial_p f(t,p) \right),\\
 &\partial_t W(t,k) = \left[ \int_{\bbR} p\,\partial_p f(t,p)\;
  \delta_{0}\bigl( \omega(k) - pk \bigr) dp \right] W(t,k),
\end{split}
\end{equation}
where $f(t,p)$ represents the background electron distribution over
momenta and $W(t,k)$ is an energy density over wave modes (here the plasma frequency,
electron mass, and speed of light are rescaled to unity). The Dirac
kernel enforces the resonance $\omega(k) = pk$, coupling a mode of
wavenumber $k$ to particles whose momentum matches its phase
velocity. Ivanov and Rudakov~\cite{ivanov1967quasilinear} take
$\omega \equiv 1$, while we take the
Bohm--Gross relation, as in Huang and Gamba~\cite{huang_gamba},
\begin{equation}\label{eq:dispersion}
  \omega(k) = (1+\lambda^{2}k^{2})^{1/2}
\end{equation}
with $\lambda \geq 0$ the rescaled Debye length, recovering Ivanov and Rudakov's choice at $\lambda = 0$. As noted by Bardos and Besse~\cite{bardos2021diffusion}, despite its 
continued use for over half a century in plasma physics~\cite{thorne_blandford,SheikhZaheerNoreenShah2023},
well-posedness of the quasilinear system remains open, even in the
one-dimensional case~\eqref{eq:system}.

\medskip\noindent
The one-dimensional system admits a formal scalar reduction,
studied by Ivanov and Rudakov~\cite{ivanov1967quasilinear}
for $\lambda=0$ and by Huang and Gamba~\cite{huang_gamba}
for the dispersion relation~\eqref{eq:dispersion}.
The resulting Cauchy--Dirichlet problem on
$\Omega=(0,\infty)$, for a nonnegative unknown $u(t,x)$
and prescribed data $u_0$ and $g$, is
\begin{equation}\label{eq:classical_form}
\begin{split}
  \diff{t} u(t,x) &= \rho_{\lambda}(x)u(t,x)\laplacian u(t,x) + \rho_{\lambda}(x)g(x)u(t,x),\ (t,x) \in (0,T) \times \Omega, \\
   u(t,0) &= 0,\ 0 < t < T, \\
   u(0,x) &= u_{0}(x),\ x \in \Omega.
\end{split}
\end{equation}
The weight $\rho_\lambda$ arises from the dispersion
relation~\eqref{eq:dispersion} through the resonance change of
variables~\eqref{eq:resonance_map} below, which relates $x$ to
the wavenumber $k$:
\begin{equation}\label{eq:rho_lambda_define}
  \rho_{\lambda}(x) \defeq (x + \lambda)\bigl((x + \lambda)^2 - \lambda^2\bigr)^{1/2}
  = x^{1/2}(x + \lambda)(x + 2\lambda)^{1/2},\ \lambda \geq 0.
\end{equation}
At $\lambda = 0$ this reads $\rho_0(x) = x^2$, for which Ivanov and
Rudakov give an asymptotic analysis of self-similar solutions. The reduction itself proceeds by inverting the
resonance condition $\omega(k) = pk$; for $\omega$ as in~\eqref{eq:dispersion}
this is the change of variables, introducing $x$ as a shifted momentum coordinate, 
\begin{equation}\label{eq:resonance_map}
  x(p) \defeq p-\lambda, \quad k(p) \defeq (p^{2}-\lambda^{2})^{-1/2},
  \quad
  p(k) \defeq (k^{-2}+\lambda^{2})^{1/2}.
\end{equation}
Without loss of generality, we consider $p \in (\lambda,\infty)$ as in~\cite{ivanov1967quasilinear,huang_gamba}:
for $\abs{p} \leq \lambda$, the resonance condition $\omega(k) = kp$ has no solutions, 
and the evolution~\eqref{eq:system} becomes trivial on this sub-resonant region, while the case $p < -\lambda$ is 
symmetric under $(p,k) \mapsto (-p,-k)$. Now~\eqref{eq:resonance_map} provides a bijection from $x = p-\lambda \in (0,\infty)$ to the wavenumber $k \in (0,\infty)$, under which the boundary $x = 0$ corresponds to
$k \to \infty$. The scalar unknown is then related to the wave spectrum $W$, particle distribution $f$,
and reaction $g$ by, 
\begin{equation}\label{eq:reconstruction}
\begin{split}
  & u(t,x) = (k(x+\lambda))^3\; W\bigl(t,\, k(x+\lambda)\bigr),\\ 
  & g(x) = \partial_x\bigl[f_0(x+\lambda)-\partial_x u_0(x)\bigr],\\
  & f(t,x+\lambda) = f_0(x+\lambda) + \partial_x(u(t,x)-u_0(x)),\quad  x > 0.
\end{split}
\end{equation}
We see that the Dirichlet condition of~\eqref{eq:classical_form} is thus a decay condition on the spectral energy $W(t,k)$ at large wavenumber. 
This derivation is given in Section~\ref{section:system}, together with a rigorous converse: a weak solution of~\eqref{eq:classical_form}
induces a weak solution of~\eqref{eq:system}.

\medskip\noindent
Problem~\eqref{eq:classical_form} lies at the intersection of two
lines of work. Nondivergence-form equations with the
prototype $\partial_t u-u\Delta u=0$ have been studied for
continuous (and smoother) data. For $\partial_t u-u\Delta u=g(u)$,
Ughi~\cite{Ughi1986} constructs nonnegative distributional
solutions in one dimension from bounded Lipschitz data; Dal Passo
and Luckhaus~\cite{dalpasso_1987} treat bounded Lipschitz domains
in higher dimensions with $u_0\in C^0(\overline{\Omega})$.
Bertsch, Dal Passo, and Ughi establish nonuniqueness phenomena
and one-sided time-derivative estimates in classical and viscosity
solution settings~\cite{BertschDalPassoUghi1990,BertschDalPassoUghi1992},
further developed by Winkler~\cite{WinklerThesis2000,WinklerJDE2003}.
More recently, Dunlap and Graham~\cite{DunlapGraham2025} obtain
uniqueness for the one-dimensional prototype in a strong-solution
class, under two-sided polynomial bounds on $u_0$ and Lipschitz
regularity of $\sqrt{u_0}$. Meanwhile, rewriting equation~\eqref{eq:classical_form}
as in~\eqref{eq:divergence_form} reveals a porous medium
diffusion term plus a gradient nonlinearity. Spatially weighted diffusion has been studied extensively
for the porous medium family \cite{KaminRosenau1981,ReyesVazquez2006,GrilloMuratoriPorzio2013},
and porous medium equations with gradient nonlinearities
have also been considered~\cite{APW}.

\medskip\noindent
We seek an existence theory under integrability assumptions,
without requiring boundedness or continuity of the initial data.
Such rough data are physically relevant: admissible initial
particle distributions $f_0$ include probability measures and atomic
data which require jump discontinuities in $u_0$ through the
reconstruction identity~\eqref{eq:reconstruction}. These data lie beyond the continuous-data theories discussed
above. Although related to weighted porous medium equations,
the nondivergence structure introduces an additional challenge
in passing to the limit in weak approximations. We therefore adapt the framework of Boccardo, Dall'Aglio,
and collaborators~\cite{boccardo_gallouet_1989,dallaglio_primary}
for quasilinear parabolic equations to our
weighted setting.

\medskip\noindent
Huang and Gamba~\cite{huang_gamba} studied existence for the scalar
problem~\eqref{eq:classical_form} from smooth, compactly
supported initial data. Their work raises three natural questions, which we address
here: existence from nonsmooth data without
compact-support restrictions, physical admissibility (e.g. nonnegativity, mass conservation) of the
reconstructed particle--wave solution, and positive-time smoothing. 
Note that nonnegativity of the scalar solution $u$ alone does
not ensure nonnegativity of the particle distribution $f$ as
the reconstruction~\eqref{eq:reconstruction} involves
$\partial_x u$.

\medskip\noindent
In this work, we construct nonnegative weak solutions
to~\eqref{eq:classical_form} on the half-line under weighted
integrability assumptions on the initial data and boundedness
of the reaction coefficient $\rho_\lambda g$
(Theorem~\ref{thm:existence}). The constructed solution becomes
locally H\"older in space and time and locally Lipschitz in space
at positive times (Theorem~\ref{thm:regularity}). This smoothing
follows from a one-sided bound on $\partial_t u$, obtained by
a time-scaling comparison in the spirit of B\'enilan--Crandall
that accommodates the spatially dependent reaction
(Proposition~\ref{prop:benilan_crandall}). Explicit solutions
with gradient jumps at the free boundary show that the spatial
regularity is sharp
(Subsection~\ref{subsection:bc_special_solutions}).
The scalar solution also induces a weak solution of the
quasilinear system~\eqref{eq:system}, preserving nonnegativity
and conserving particle mass under natural hypotheses on the
initial data (Theorem~\ref{thm:system}). The induced solution
inherits positive-time regularity: admissible initial particle
measures become locally bounded densities, and the wave spectrum
decays quantitatively at large wavenumber.

\medskip\noindent
We do not address uniqueness, which we expect to fail at least for the 
scalar problem (Remark~\ref{rmk:uniqueness}). Physical admissibility is established for the particular
solution we construct. Precise hypotheses and statements
are given in Subsection~\ref{subsection:main_results}.

\subsection{Notation}\label{subsection:notation}
To avoid unnecessary bookkeeping of real-valued sequences which are asymptotically zero, we denote such quantities by an equivalence class $\omega(n)$, writing 
\begin{equation}\label{eq:omega_define} a_n \equiv \omega(n) \iff \limsup_{n}\ \abs{a_n} = 0. \end{equation}
We write the spatial domain $\Omega \defeq (0,\infty)$, the parabolic cylinder $Q_T \defeq (0,T) \times \Omega$, parabolic boundary $\partial Q_{T} \defeq \overline{Q_T} \setminus(0,T] \times \Omega$, 
and for $t_0 > 0$ we write $Q_{(t_0,T)} \defeq (t_0,T)\times \Omega$. We define the Borel measure $\nu_{\lambda}$ on $(0,\infty)$ by its Lebesgue density,
\begin{equation*}d\nu_{\lambda} \defeq \frac{1}{\rho_{\lambda}(x)}dx = \begin{cases} \frac{1}{x^{1/2}(x+\lambda)(x+2\lambda)^{1/2}}dx & \lambda > 0\\ \frac{1}{x^2}dx & \lambda = 0. \end{cases} \end{equation*}
Any integral is with respect to Lebesgue measure on its domain unless otherwise indicated, i.e.,
$\fintegral{Q_T} f \defeq \fintegral{Q_T} f(x,t)dtdx$. $L^{p}(A)$ denotes the usual Lebesgue space for a measurable subset $A \subset \bbrd$, 
and $L^{p}(A;d\mu)$ the corresponding space with measure $\mu$. For weighted Sobolev spaces,
following~\cite{muratori_primary}, we define the norm 
\begin{equation*}
\norm{f}_{V^p(a,b;d\nu_\lambda)}
\defeq \norm{f}_{L^p((a,b);d\nu_\lambda)}
+\norm{\partial_x f}_{L^2(a,b)}.
\end{equation*}
Write $\mathcal{D}(a, b) := C^\infty_c(a, b)$ and
\[
  \mathcal{D}_L(a, b) := \{f \in C^\infty((a, b)) \cap C^{0}([a,b]) : \lim_{x \to a^{+}} f(x) = 0\}
\]
to define
\begin{equation*}V^p_0((a, b)) := \overline{\mathcal{D}(a, b)}^{\norm{\cdot}_{V^p(a,b; d\nu_\lambda)}}, 
  \quad 
  V^p_L((a, b)) := \overline{\mathcal{D}_L(a, b)}^{\norm{\cdot}_{V^p(a,b; d\nu_\lambda)}}.
\end{equation*}
We suppress the dependence on $\nu_\lambda$ in the notation
for the resulting spaces $V_0^p$ and $V_L^p$.
Note that when $b = \infty$ one has $V^{p}_L((0,\infty)) = V^{p}_0((0,\infty))$. We record useful embeddings of these spaces and their associated Hardy inequalities 
in Appendix~\ref{appendix:hardy_sobolev}. Note that for all $\lambda \geq 0$, $1 \leq p < \infty$, and $0 < b < \infty$, monotonicity of $\rho_\lambda$ implies
\begin{equation}\label{eq:lambda_to_lebesgue}
  \norm{f}_{L^p(0,b)} \leq \rho_\lambda(b)^{1/p} \norm{f}_{L^p(0,b; d\nu_\lambda)}, 
  \quad 
  \norm{f}_{L^p(b,\infty; d\nu_\lambda)} \leq \rho_\lambda(b)^{-1/p} \norm{f}_{L^p(b,\infty)}.
\end{equation}
We follow Krylov's notation~\cite{krylov} for parabolic H\"older spaces: for $\alpha \in (0,1)$, 
$C^{\alpha/2,\alpha}(\mc{O})$ on a relatively compact $\mc{O} \Subset Q_T$ is defined with the seminorm 
\[ [f]_{C^{\alpha/2,\alpha}(\mc{O})} := \sup_{\substack{(t_1,x_1),\, (t_2,x_2) \in \mc{O} \\ (t_1,x_1) \neq (t_2,x_2)}}
\frac{\abs{f(t_1,x_1)-f(t_2,x_2)}}{\abs{x_1-x_2}^{\alpha}+\abs{t_2-t_1}^{\alpha/2}}, \]
with $C^{k+\alpha/2,2k+\alpha}$ for $k \in \bbN$ defined in the usual way. For $O \Subset \Omega$, 
$C^{k+\beta}(O)$ are standard H\"older spaces with $k$-th derivative having finite $\beta$-seminorm.

\subsection{Main Results}\label{subsection:main_results}
We now present our weak-solution framework for~\eqref{eq:classical_form}. First, we transfer the weight $\rho_\lambda$ to the
left-hand side of~\eqref{eq:classical_form}, incorporating $\rho_\lambda^{-1}$ into appropriate function spaces through the
measure $\nu_\lambda$, analogous to Muratori's method for weighted porous
medium equations~\cite{muratori_primary}. Then, we rewrite the
leading-order term in divergence form up to a gradient nonlinearity, using 
$u\,\partial_x^2 u = \partial_x^2\bigl(\tfrac{u^2}{2}\bigr)
- \abssq{\partial_x u}$. For smooth solutions, \eqref{eq:classical_form} becomes
\begin{equation}\label{eq:divergence_form}
  \rho_\lambda^{-1}(x)\,\partial_t u
  = \partial_x^2\Bigl(\frac{u^2}{2}\Bigr) - \abssq{\partial_x u} + gu.
\end{equation}
Testing \eqref{eq:divergence_form} and integrating by parts motivates the following. 
\begin{definition}\label{defn:weak_solution}
We say a nonnegative function $u \colon Q_T \longrightarrow \bbR$ is a
\textbf{weak solution} to~\eqref{eq:classical_form} with nonnegative
initial data $u_0\in L^{1}_{loc}(\Omega;d\nu_{\lambda})$ if
$u^2 \in L^{3/2}(0,T;V_0^{3/2}(\Omega))$,
$\partial_{x}u^2 \in L^2(Q_T)$,
$\partial_x u\in L^2((0,T')\times D)$ for every
$0<T'<T$ and every $D\Subset\Omega$, and for all
$\phi \in C^{\infty}_c([0,T) \times \Omega)$, $u$ satisfies
\begin{equation}\label{eq:weak_form}
  \begin{split}
    & -\fintegral{Q_T} u(t,x)\diff{t}\phi(t,x)\,dt\,d\nu_{\lambda}
    + \fintegral{Q_T}\partial_{x}{\left(\frac{u^2}{2}\right)}
      \cdot\partial_{x}{\phi}\,dt\,dx
    + \fintegral{Q_T}\abs{\partial_{x}{u}}^2\phi \,dt\,dx\\
    & = \fintegral{\Omega} u_0(x)\phi(0,x)\,d\nu_{\lambda}
    + \fintegral{Q_T} \phi(t,x)u(t,x)g(x)\,dt\,dx.
  \end{split}
\end{equation}
\end{definition}

All terms in~\eqref{eq:weak_form} are finite: the weighted $L^3$
bound controls the time and reaction terms, the stated gradient
bounds control the spatial terms on the support of $\phi$, and
the initial datum is locally integrable. The boundary condition
is incorporated through $u^2\in L^{3/2}(0,T;V_0^{3/2}(\Omega))$.

\begin{theorem}[Existence of Weak Solutions]\label{thm:existence}
For any $\lambda \geq 0$ and any $T > 0$, suppose
$\rho_{\lambda}g \in L^{\infty}(0, \infty)$, and suppose $u_0 \geq 0$
satisfies
\begin{equation*}u_0 \in
\begin{cases}
L^{\delta}(0, \infty;d\nu_{0}) \cap L^{3}(0, \infty;d\nu_{0}),
  & \text{for some } 0 < \delta < 1, \quad \text{if } \lambda = 0 \\
L^{1}(0, \infty;d\nu_{\lambda}) \cap L^3(0, \infty;d\nu_{\lambda}),
  & \text{if } \lambda > 0.
\end{cases}
\end{equation*}
Then there exists a nonnegative weak solution $u$
of~\eqref{eq:classical_form} on $Q_T$ with initial data $u_0$, in the
sense of Definition~\ref{defn:weak_solution}, satisfying
$\partial_{x}u \in L^2(Q_T)$, and
\[
  u \in L^{\infty}(0,T;L^{p}(\Omega;d\nu_{\lambda}))\ \forall\
  p \in [1,3],
\]
with
$u \in L^{\infty}(0,T;L^{\delta}(\Omega;d\nu_{\lambda}))$ when
$\lambda = 0$.
\end{theorem}

The requirement $u_0 \in L^{\delta}(d\nu_0)$ when $\lambda = 0$, employed in the existence argument
(Propositions~\ref{prop:evolution_compactness}
and~\ref{prop:gradient_tail_bounds}), is a
local condition near $x = 0$: since
$\nu_0((b,\infty)) < \infty$ for $0 < b < \infty$, one has
\[
  u_0 \in L^\delta((0,b); d\nu_0) \cap L^3((0,\infty); d\nu_0) \implies u_0 \in L^{\delta}(d\nu_0).
\]
Note also that $u_0 \in L^{\delta}(d\nu_0)$ is a strictly weaker assumption than imposing a higher polynomial weight at the
origin, as
\[
  L^{1}\bigl((0,b); x^{-2-\epsilon}dx\bigr) \hookrightarrow
  L^{\delta}\bigl((0,b); x^{-2}dx\bigr)
\]
for $\delta = 2/(2+\epsilon)$, $\epsilon \in (0,1)$. When $\lambda > 0$,
$\nu_\lambda$ is finite and the condition is implied by
$u_0 \in L^1(d\nu_\lambda)$ for every $\delta \in (0,1)$, a freedom we
use in the sequel.

\begin{remark}[Nonuniqueness]\label{rmk:uniqueness}
We construct a particular weak solution and do not address
uniqueness for~\eqref{eq:classical_form}, which we expect to
fail in general within our weak-solution class.
For the homogeneous prototype $\partial_t u=u\,\partial_x^2u$,
nonuniqueness is known for suitable initial data, both through
shrinking positivity sets and through different waiting times
for positivity at isolated zeros
\cite{dalpasso_1987,Ughi1986,BertschDalPassoUghi1992}.
For our constructed solution, the bound~\eqref{eq:bc_for_u}
and positive-time continuity imply that the positivity set
is nondecreasing on $(0,T]$, excluding shrinking positivity
sets and extinction after a positive time.
This does not settle uniqueness: the prototype admits distinct
solutions satisfying the analogous one-sided bound
\cite[Theorem~3.1]{BertschDalPassoUghi1992}.
We do not establish nonuniqueness for the weighted equation
considered here.
\end{remark}

We prove that the solution of Theorem~\ref{thm:existence} gains regularity at positive times.
\begin{theorem}[Positive-Time Smoothing]\label{thm:regularity}
The weak solution $u$ from Theorem~\ref{thm:existence} admits a
representative that is locally H\"older continuous in
$(0, T] \times \Omega$ and locally Lipschitz in $x$: for every
$0 < t_0 < T$ and every $D \Subset \Omega$,
\begin{equation*}
  \norm{u}_{C^{1/4,\, 1/2}((t_0, T) \times D)} < \infty
  \quad \text{and} \quad
  \norm{\partial_x u}_{L^\infty((t_0, T) \times D)} < \infty.
\end{equation*}
Quantitative bounds on these norms in terms of the data are given in
the proof in Section~\ref{section:regularity}.
\end{theorem}
The smoothing mechanism driving Theorem~\ref{thm:regularity} is a
one-sided lower bound on $\partial_t u$ of B\'enilan--Crandall type. We use this name 
as the estimate follows from time-scaling and comparison in the spirit of~\cite{BenilanCrandall1981},
and the divergence form~\eqref{eq:divergence_form} expresses the equation as a porous
medium term plus a quadratic-gradient term; such estimates are classical for the 
porous medium equation and its inhomogeneous analogues \cite{vazquez,muratori_secondary}. Analogous bounds for
the homogeneous nondivergence-form prototype are known
from~\cite{BertschDalPassoUghi1992,WinklerJDE2003}; here we establish
the estimate for the constructed weak solution with the spatial
weight and reaction. The negative part of the reaction
coefficient enters through the parameter $\theta$ and the associated
time scale $K_\theta(t)$: 
\begin{align}
& \label{eq:defn_theta}
\theta \defeq \max\bigl(\esssup_{\Omega} (-\rho_{\lambda} g),\ 0\bigr), \\
& \label{eq:defn_k_theta}
 K_\theta(t) \defeq \begin{cases}
    t, & \theta = 0, \\
    \dfrac{1}{\theta}\dfrac{\exp(\theta t)-1}{\exp(\theta t)},
      & \theta > 0.
                  \end{cases}
\end{align}
Since $\rho_\lambda g \in L^\infty$, $\theta \in [0,\infty)$ and
$K_\theta(t) > 0$ for $t > 0$. 
\begin{proposition}[B\'enilan--Crandall Bound]\label{prop:benilan_crandall}
With $\theta$ and $K_{\theta}$ as in~\eqref{eq:defn_theta}
and~\eqref{eq:defn_k_theta}, the weak solution $u$ from
Theorem~\ref{thm:existence} satisfies, in the sense of distributions
in $Q_T$,
\begin{equation}\label{eq:bc_for_u}
  \partial_t u \geq -\frac{u}{K_\theta(t)}.
\end{equation}
\end{proposition}
The spatial Lipschitz regularity in
Theorem~\ref{thm:regularity} is sharp within the weak-solution
class of Definition~\ref{defn:weak_solution}:
the explicit solutions of
Subsection~\ref{subsection:bc_special_solutions} have gradients
that jump to zero at the free boundary and attain equality
in~\eqref{eq:bc_for_u} for both $\theta=0$ and $\theta>0$.

\medskip\noindent
Finally, we construct a weak solution for the quasilinear system \eqref{eq:system} from the scalar
solution of Theorem~\ref{thm:existence}. Define
\begin{equation}\label{eq:coupling_profile}
  h(p) \defeq -\int_p^\infty g(p'-\lambda)\,dp',
  \qquad p \in (\lambda,\infty),
\end{equation}
which is well defined under the boundedness assumption on
$\rho_\lambda g$. With $x(k)=p(k)-\lambda$ from~\eqref{eq:resonance_map}, we define the induced particle-wave pair 
following \eqref{eq:reconstruction}, 
\begin{equation}\label{defn:induced_pair}
  W^\ast(t,k) \defeq k^{-3}u(t,x(k)),
  \qquad
  f^\ast(t,p) \defeq h(p)+\partial_x u(t,p-\lambda),
\end{equation}
where the derivative is understood weakly. The corresponding
initial data are
$W^\ast_0(k)=k^{-3}u_0(x(k))$ and
$f^\ast_0=h+\partial_p[u_0(p-\lambda)]$, the latter defined
in $\mathcal{D}'((\lambda,\infty))$.

While $W^\ast\geq0$ follows from $u\geq0$, preservation of
particle nonnegativity requires
\begin{equation}\label{eq:f0_positive_hypothesis}
  f^\ast_0 \geq 0
  \quad\text{in }\mathcal{D}'((\lambda,\infty)).
\end{equation}
This allows $f^\ast_0$ to be a locally finite measure;
in particular, $u_0\in BV_{\mathrm{loc}}(\Omega)$ follows
from this condition rather than being assumed separately.
To establish conservation of particle mass, we require the additional assumptions 
\begin{equation}\label{eq:moment_hypotheses}
  \int_0^\infty x\,|g(x)|\,dx<\infty,
  \qquad
  \int_0^\infty x\,u_0(x)^3\,d\nu_\lambda(x)<\infty,
\end{equation}
and normalize
\begin{equation}\label{eq:unit_mass}
  \int_\lambda^\infty h(p)\,dp=1.
\end{equation}
The unit normalization is a convention; the same argument
conserves any prescribed positive mass.

\medskip\noindent
The weak formulation of the system~\eqref{eq:system} is as follows
(cf.~\cite[Eqs.~(3.1), (3.3)]{huang_gamba}): for
$\zeta\in C^\infty_c((\lambda,\infty)\times[0,T))$ and
$\xi\in C^\infty_c((0,\infty)\times[0,T))$,
\begin{equation}\label{eq:system_weak_intro}
\begin{split}
  &-\int_0^T\!\!\int_\lambda^\infty
    f^\ast\,\partial_t\zeta\,dp\,dt
  -\bigl\langle f^\ast_0,\zeta(\cdot,0)\bigr\rangle
  =\int_0^T\!\!\int_\lambda^\infty
    f^\ast\,\partial_p\bigl(\gamma u\,\partial_p\zeta\bigr)
    \,dp\,dt,\\
  &-\int_0^T\!\!\int_0^\infty
    W^\ast\,\partial_t\xi\,dk\,dt
  -\int_0^\infty W^\ast_0\,\xi(k,0)\,dk\\
  &\qquad=-\int_0^T\!\!\int_\lambda^\infty
    f^\ast\,\partial_p\bigl(p\gamma u\,\xi(k(p),t)\bigr)
    \,dp\,dt,
\end{split}
\end{equation}
where $\gamma(p)\defeq p(p^2-\lambda^2)^{1/2}$ and $u$
is evaluated at $(t,p-\lambda)$. The initial particle term
is understood as a distributional pairing.

\begin{theorem}[Induced Solution of the System]\label{thm:system}
Let $u_0$ and $g$ satisfy the hypotheses of
Theorem~\ref{thm:existence}. Then:
\begin{enumerate}
  \item[(i)] the induced pair $(f^\ast,W^\ast)$
  of~\eqref{defn:induced_pair} is a weak solution of
  the system~\eqref{eq:system} in the sense
  of~\eqref{eq:system_weak_intro}, with $W^\ast\geq0$;
  \item[(ii)] if moreover~\eqref{eq:f0_positive_hypothesis}
  holds, then $u$ may be constructed so that additionally
  $f^\ast(t,\cdot)\geq0$ a.e.\ on $(\lambda,\infty)$
  for a.e.\ $t\in(0,T)$;
  \item[(iii)] if additionally~\eqref{eq:moment_hypotheses}
  and~\eqref{eq:unit_mass} hold, the same construction
  conserves mass,
  $\int_\lambda^\infty f^\ast(t,p)\,dp=1$
  for a.e.\ $t\in(0,T)$.
\end{enumerate}
\end{theorem}

Parts (ii) and (iii) rely on our choice for regularized
data and are not asserted for every scalar weak solution.
The induced pair also inherits positive-time regularity:
admissible initial particle measures become locally bounded
densities, while $W^\ast$ gains local H\"older and spatial
Lipschitz regularity, persistence of its excited spectrum,
and quantitative decay at large wavenumber
(Proposition~\ref{prop:regularity_induced_pair} and
Corollary~\ref{cor:positive_time_properties}).
Subsection~\ref{subsection:system_examples} illustrates these
results with a bump-on-tail configuration whose initial particle
distribution consists of two Dirac masses and an explicit
separable solution relaxing at the rate $\theta$.

\subsection{Strategy of Proof}\label{subsection:strategy}

We prove Theorem~\ref{thm:existence} by passing to the limit
in the regularized problems of
Section~\ref{section:regularized_problems}. An Aubin--Lions
lemma adapted to the weighted framework of
Definition~\ref{defn:weak_solution} provides compactness
of the regularized solutions
(Subsection~\ref{subsection:compactness}). The principal difficulty is passage to the limit in the
quadratic-gradient term. We adapt the approach of Boccardo,
Dall'Aglio, and collaborators
\cite{boccardo_1989,boccardo_2002,bdgo1997,dallaglio_primary}
to degeneracy at $u=0$, the spatial weight $\rho_\lambda$,
and the unbounded half-line. The key step is convergence of
the truncated gradients
(Proposition~\ref{prop:trunc_gradients_converge}), leading
to strong local gradient convergence
(Corollary~\ref{cor:gradients_converge_locally}).
The technical calculations are collected in
Appendix~\ref{appendix:dall_aglio}.

\medskip\noindent
The B\'enilan--Crandall bound follows from time-scaling
and comparison for the regularized problems and passes
to the limit, yielding the positive-time regularity
of Theorem~\ref{thm:regularity}
(Section~\ref{section:regularity}).
Testing the scalar weak formulation with test functions
pulled back through the resonance map gives the weak
formulation (Section~\ref{section:system}).
We prove that a suitable choice of regularized data ensures nonnegativity
of the induced particle distribution
(Proposition~\ref{prop:f_nonneg}).

\section{Regularized Problems}\label{section:regularized_problems}
We construct a family of solutions $\{u_n\}_{n \in \bbN}$ to regularized problems on bounded domains and identify a candidate weak solution $u$ to~\eqref{eq:classical_form} as a subsequential limit of $\{u_n\}_{n \in \bbN}$. 
Subsection~\ref{subsection:uniform_estimates} establishes a priori bounds (uniform in $n$) while Subsection~\ref{subsection:compactness} extracts a limit by compactness.
\subsection{Uniform Estimates}\label{subsection:uniform_estimates}
\begin{definition}[Regularized Solutions $u_n$]\label{defn:n_problem}
Fix $\lambda \geq 0$. For each $n \in \bbN$, choose a domain $\Omega_n = (a_n, b_n)$
satisfying
\begin{equation}\label{eq:omega_n_requirement}
\begin{gathered}
  a_{n+1} < a_n, \quad b_n < b_{n+1}, \\
  \lim_{n \to \infty} a_n = 0, \quad \lim_{n \to \infty} b_n = \infty, \quad
  \lim_{n \to \infty} n^{-\delta} \nu_\lambda(\Omega_n) = 0.
\end{gathered}
\end{equation}
($(a_n, b_n) = (n^{-\delta/2}, n)$ works for all $\lambda \geq 0$ for large $n$), and write $Q_{T,n} \defeq (0,T) \times \Omega_{n}$ 
as the $n$-parabolic cylinder. Choose data $u_{0, n}$ and $g_n$ 
approximating $u_0$ and $g$, satisfying 
\begin{equation}\label{eq:n_problem_data_requirement}
\begin{split}
  &\rho_{\lambda} g_n \in C^\infty_c(\Omega_n), u_{0,n} \geq 0, \quad u_{0, n} \in C^{2+\beta}_C(\Omega_n) 
  \text{ for some } \beta \in (0, 1), \\
  &u_{0, n} \to u_0 \text{ in } L^p(\Omega; d\nu_\lambda) \text{ for } p \in [\delta, 3], \\
  &\rho_{\lambda} g_n \to \rho_{\lambda} g \text{ in } L^p_{loc}(\Omega) \text{ and Lebesgue-a.e.}, \\
  &\rho_{\lambda} g_n \geq -\theta, \quad
  \norm{g_n \rho_{\lambda}}_{L^\infty(\Omega)} \leq \norm{g \rho_{\lambda}}_{L^\infty(\Omega)}.\\
  &\support u_{0,n} \cup \support(\rho_\lambda g_n) \subset [a_n + r_n,\ b_n - r_n],
\end{split}
\end{equation}
for $r_n > 0$ satisfying $r_n \to 0^+$ and $r_n n^\delta \to \infty$ 
($r_{n} = n^{-\delta/2}$ suffices), and where $\delta$ is from Theorem~\ref{thm:existence} (see the discussion following its statement).
Here $\theta$ is as defined in~\eqref{eq:defn_theta}. Define $u_n$ as the unique classical solution to
\begin{equation}\label{eq:n_problem}
\begin{split}
  \diff{t} u_n &= \rho_{\lambda} (u_n + \tfrac{1}{n}) \diff{x}^2 u_n 
    + \rho_{\lambda} g_n (u_n + \tfrac{1}{n}) \quad \text{on } (0, T) \times \Omega_n, \\
  u_n(0, x) &= u_{0, n}(x), \quad u_n(t, a_n) = u_n(t, b_n) = 0.
\end{split}
\end{equation}
\end{definition}
Define also
\begin{equation}\label{eq:defn_tilde_un}
  \tilde{u}_n \defeq u_n + \tfrac{1}{n} \quad \text{on } \overline{Q_{T,n}}.
\end{equation}
We verify existence and uniqueness of classical solutions $u_n$ by standard quasilinear theory~\cite{Ladyenskaja1968}.
\begin{proposition}\label{prop:n_are_classical}
For each $n$, there exists a unique classical solution
\[
  u_n \in C^{1+\beta/2, 2+\beta}(\overline{Q_{T,n}})
\]
to~\eqref{eq:n_problem}. The extension by zero $u_n$ to $Q_T$, not relabeled, 
is continuous and locally Lipschitz on $\overline{Q_T}$, satisfying 
\begin{equation}\label{eq:n_maximum_principle}
  -\frac{1}{n}+\frac{e^{-\theta t}}{n} \leq u_n(t, x)
  \leq \left(\norm{u_{0,n}}_{L^\infty(\Omega_n)}+\frac1n\right)
      e^{t\linfty{g\rho_\lambda}}-\frac1n.
\end{equation}
In particular,
\begin{equation}\label{eq:shifted_positivity}
  \tilde{u}_n(t,x)\geq \frac{e^{-\theta t}}n>0
  \qquad\text{on }\overline{Q_{T,n}}.
\end{equation}
With $r_n \to 0$ as in Definition~\ref{defn:n_problem},
\begin{equation}\label{eq:n_gradient_vanishes}
  \sup_{t \in [0, T]} \bigl( \abs{\diff{x} u_n(t, a_n)} + \abs{\diff{x} u_n(t, b_n)} \bigr) 
  \leq \frac{2(\norm{u_{0,n}}_{L^\infty(\Omega_n)}+1/n)e^{T\linfty{g\rho_\lambda}}}{r_n}.
\end{equation}
The one-sided bounds, uniform in the size of $u_{0,n}$, are
\begin{equation}\label{eq:n_boundary_one_sided}
  \diff{x}u_n(t,a_n)\geq-\frac1{nr_n},
  \qquad \diff{x}u_n(t,b_n)\leq\frac1{nr_n}.
\end{equation}
\end{proposition}
\begin{proof}
The shifted function satisfies
\[
  \partial_t \tilde{u}_n=\rho_\lambda \tilde{u}_n\partial_x^2\tilde{u}_n
       +\rho_\lambda g_n \tilde{u}_n,\qquad
  \tilde{u}_n(0,\cdot)=u_{0,n}+1/n,\qquad \tilde{u}_n|_{\partial\Omega_n}=1/n.
\]
Put $M=\linfty{g\rho_\lambda}$. The spatially constant functions
\[
  n^{-1}e^{-\theta t} \quad\text{and}\quad
  (\norm{u_{0,n}}_{L^\infty(\Omega_n)}+1/n)e^{Mt}
\]
are respectively a
subsolution and a supersolution, including on the parabolic boundary.
Comparison~\cite[Theorem I.2.9]{Ladyenskaja1968} gives
\eqref{eq:n_maximum_principle}--\eqref{eq:shifted_positivity}.
For existence, apply the quasilinear theory of
\cite[Theorems V.6.1 and IV.5.2]{Ladyenskaja1968} first with a smooth,
uniformly positive extension of the diffusion coefficient below
$\tilde{u}_n=e^{-\theta T}/(2n)$. The comparison bounds keep the solution in a
bounded interval where this extension agrees with $\rho_\lambda \tilde{u}_n$;
since $\min_{\overline{\Omega_n}}\rho_\lambda>0$, the problem is uniformly
parabolic there. The cited theory and continuation give the unique
classical solution on $[0,T]$. The compatibility condition at
$\partial\Omega_n$ follows from~\eqref{eq:n_problem_data_requirement},
and the zero boundary values give the locally Lipschitz extension.

For the boundary estimates put
$U_{n,T}=(\norm{u_{0,n}}_{L^\infty(\Omega_n)}+1/n)e^{MT}$.
On the left strip $[a_n,a_n+r_n]$, where both $u_{0,n}$ and $g_n$
vanish, the stationary affine barriers
$-(x-a_n)/(nr_n)$ and $U_{n,T}(x-a_n)/r_n$ bound $u_n$ from below
and above on the parabolic boundary. Comparison therefore gives
\[
  -\frac1{nr_n}\leq \partial_xu_n(t,a_n)\leq\frac{U_{n,T}}{r_n}.
\]
The corresponding barriers on $[b_n-r_n,b_n]$ give
$-U_{n,T}/r_n\leq\partial_xu_n(t,b_n)\leq1/(nr_n)$.
Since $U_{n,T}\geq1/n$, these imply
\eqref{eq:n_gradient_vanishes} and~\eqref{eq:n_boundary_one_sided}.
\end{proof}
We record a local-in-time weak form of~\eqref{eq:n_problem} to estimate solutions $u_n$.
\begin{lemma} For any $t \geq 0$ and $\psi \in C^{\infty}_c([0,T) \times \Omega_n)$
\begin{equation} \label{eq:energy_form_n} 
  \begin{split}
   & \fintegral{\Omega_n} \diff{t}\frac{u_{n}}{\rho_{\lambda}(x)} \psi(t,x)
   + \fintegral{\Omega_n} \left(u_{n}+\frac{1}{n}\right)\partial_{x}{u_{n}}\partial_{x}{\psi(t,x)} + \fintegral{\Omega_n} \psi(t,x)\abssq{\partial_{x}{u_{n}}}\\ 
   &  = \fintegral{\Omega_n} \left(u_n+\frac{1}{n}\right)g_{n}\psi(t,x).
  \end{split} 
\end{equation} 
\end{lemma}
\begin{proof}
  Divide~\eqref{eq:n_problem} through by $\rho_{\lambda}$, test with $\psi(t,x)$, and integrate over $\Omega_n$. 
\end{proof} 
\begin{proposition}[Uniform Estimates]\label{prop:uniform_estimates}
For any $m \in \{0, 1, 2\}$, $T > 0$, $\beta \defeq 1-\delta$ with $\delta$ from Theorem~\ref{thm:existence} (and any $\delta \in (0,1)$ when $\lambda > 0$)
there exist constants $C_{m+1}$, $C_\delta$, and $C'_T$ (depending only on $g,T$ and $u_0$) such that the family 
of solutions $\{u_n\}_{n \in \mathbb{N}}$ from~\eqref{eq:n_problem} satisfies, uniformly in $n \in \mathbb{N}$,
\begin{align}
& \sup_{0 \leq t \leq T} \norm{u_n(\cdot, t)}_{L^{m+1}(\Omega; d\nu_\lambda)}^{m+1}
    + \norm{\partial_{x}{} \abs{u_n}^{m/2 + 1}}_{L^2(Q_T)}^2 
    \leq C_{m+1} + \omega(n), \label{eq:lp_apriori} \\[4pt]
& \sup_{0 \leq t \leq T} \norm{u_n + \tfrac{1}{n}}_{L^{\delta}(\Omega_n; d\nu_\lambda)}^\delta
    + \tfrac{\delta^2}{(1 - \beta/2)^2} \fintegral{Q_{T,n}}
    \abssq{\partial_{x}{}(u_n + \tfrac{1}{n})^{1 - \beta/2}}
    \leq C_\delta + \omega(n), \label{eq:delta_apriori} \\[4pt]
&  4 \norm{t \diff{t}(u_n + \tfrac{1}{n})^{1/2}}_{L^2(0, T; L^{2}(\Omega; d\nu_\lambda))}^2
    + \tfrac{T}{2} \fintegral{\Omega} \normsq{\partial_{x}{u_n}(T, \cdot)}
    \leq C'_T + \omega(n). \label{eq:l2_time_apriori}
\end{align}
Explicit values of these constants are given in the proof below.
\end{proposition} 
\begin{proof}
Put $M=\linfty{g\rho_\lambda}$ and, for $p\in\{\delta,1,2,3\}$, define
\[
 E_{p,n}(t)=\int_{\Omega_n}\tilde{u}_n(t,x)^p\,d\nu_\lambda(x).
\]
Testing~\eqref{eq:n_problem}, divided by $\rho_\lambda$, with
$p \tilde{u}_n^{p-1}$ and retaining the boundary term gives the identity
\begin{equation}\label{eq:shifted_power_energy}
  \frac{d}{dt}E_{p,n}
  +p^2\int_{\Omega_n}\tilde{u}_n^{p-1}\abssq{\partial_xu_n}
  =p\int_{\Omega_n}g_n \tilde{u}_n^p
   +\frac{p}{n^p}[\partial_xu_n]_{a_n}^{b_n}.
\end{equation}
All these tests are legitimate by~\eqref{eq:shifted_positivity} and
classical regularity. By~\eqref{eq:n_boundary_one_sided}, the last term
is bounded above by $2p/(n^{p+1}r_n)=\omega(n)$, uniformly in time.
Thus Gronwall and integration on $[0,T]$ give
\begin{equation}\label{eq:shifted_power_bound}
  \sup_{0\leq t\leq T}E_{p,n}(t)
  +p^2\int_{Q_{T,n}}\tilde{u}_n^{p-1}\abssq{\partial_xu_n}
  \leq 2e^{pMT}E_{p,n}(0)+\omega(n).
\end{equation}
The domain condition~\eqref{eq:omega_n_requirement}
implies $n^{-p}\nu_\lambda(\Omega_n)\to0$ for all $p\geq\delta$.
Together with~\eqref{eq:n_problem_data_requirement}, Minkowski's
inequality for $p\geq1$ and subadditivity of $s\mapsto s^\delta$ give
\[
  E_{p,n}(0)=\norm{u_0}_{L^p(\Omega;d\nu_\lambda)}^p+\omega(n).
\]

\textit{Proof of~\eqref{eq:lp_apriori}:}
Take $p=m+1\in\{1,2,3\}$ in~\eqref{eq:shifted_power_bound} and use
\[
 \tilde{u}_n^{p-1}\abssq{\partial_xu_n}
 =\frac{4}{(p+1)^2}\abssq{\partial_x \tilde{u}_n^{(p+1)/2}},
 \qquad \frac{4p^2}{(p+1)^2}\geq1.
\]
Since $\abs{u_n}^p\leq \tilde{u}_n^p+n^{-p}$ on $\Omega_n$, the same bound
controls the required $L^p$ norm of $u_n$. For the gradient term,
$\abssq{\partial_x\abs{u_n}}=\abssq{\partial_xu_n}$ a.e.; for
$m=1,2$ use respectively
$\abs{u_n}\leq \tilde{u}_n+1/n$ and $u_n^2\leq \tilde{u}_n^2+1/n^2$.
The additional gradient terms are $\omega(n)$ by the $p=1$ bound.
Consequently~\eqref{eq:lp_apriori} holds with
\[
  C_p=2e^{pMT}\norm{u_0}_{L^p(\Omega;d\nu_\lambda)}^p,
  \qquad p\in\{1,2,3\}.
\]
For use in compactness, the $p=1$ identity also controls the boundary
fluxes in absolute value. Indeed,~\eqref{eq:n_boundary_one_sided}
and~\eqref{eq:shifted_power_energy} give
\begin{equation}\label{eq:n_boundary_flux_l1}
\begin{split}
 &\frac1n\int_0^T\bigl(\abs{\partial_xu_n(t,a_n)}
                 +\abs{\partial_xu_n(t,b_n)}\bigr)\,dt\\
 &\qquad\leq E_{1,n}(0)+M\int_0^T E_{1,n}(t)\,dt
                   +\frac{4T}{n^2r_n}\leq C+\omega(n).
\end{split}
\end{equation}
This uses no uniform $L^\infty$ bound on the initial approximations.

\textit{Proof of~\eqref{eq:delta_apriori}:}
Take $p=\delta$ in~\eqref{eq:shifted_power_bound}. Since
$\beta=1-\delta$,
\[
  \tilde{u}_n^{-\beta}\abssq{\partial_xu_n}
  =(1-\beta/2)^{-2}
       \abssq{\partial_x \tilde{u}_n^{1-\beta/2}},
\]
which gives~\eqref{eq:delta_apriori} with
$C_\delta=2e^{\delta MT}\norm{u_0}_{L^\delta(\Omega;d\nu_\lambda)}^\delta$.
For $\lambda>0$, the finiteness of $\nu_\lambda(\Omega)$ ensures that
this norm is finite under the hypotheses of Theorem~\ref{thm:existence}.

\textit{Proof of~\eqref{eq:l2_time_apriori}:}
For $\xi\in C^1([0,T])$, $\xi\geq0$, test~\eqref{eq:energy_form_n}
with $\psi=\xi(t)\partial_tu_n/\tilde{u}_n$. This test vanishes at
$\partial\Omega_n$ because $\partial_tu_n=0$ there.
The terms containing $\abssq{\partial_xu_n}$ cancel, and integration
by parts in time and space gives
\[
\begin{split}
 &4\int_{Q_{T,n}}\frac{\xi(t)}{\rho_\lambda}
          \abssq{\partial_t \tilde{u}_n^{1/2}}
       +\frac{\xi(T)}2\int_{\Omega_n}\abssq{\partial_xu_n(T,x)}\\
 &=\frac{\xi(0)}2\int_{\Omega_n}\abssq{\partial_xu_n(0,x)}
       +\frac12\int_{Q_{T,n}}\xi'(t)\abssq{\partial_xu_n}\\
 &\quad+\left[\xi(t)\int_{\Omega_n}g_nu_n(t,x)\,dx\right]_0^T
       -\int_{Q_{T,n}}\xi'(t)g_nu_n.
\end{split}
\]
Taking $\xi(t)=t$ eliminates the initial-gradient term. By
\eqref{eq:lp_apriori}, the right-hand side is at most
$(\tfrac12+2TM)C_1+\omega(n)$, using
$\abs{\int g_nu_n}\leq M\norm{u_n}_{L^1(d\nu_\lambda)}$.
Since $t^2\leq Tt$ on $[0,T]$, this implies
\eqref{eq:l2_time_apriori} with
\[
 C'_T=\max\{1,T\}\bigl(\tfrac12+2TM\bigr)C_1.
\]
\end{proof}

\smallskip\noindent\textit{Control of the negative part.}
Write $(u_n)_-\defeq\max\{-u_n,0\}$. The estimates above imply
\begin{equation}\label{eq:negative_part_estimates}
\begin{split}
 &\norm{(u_n)_-}_{L^\infty(Q_T)}\leq\frac1n,\\
 &\norm{\partial_x(u_n)_-}_{L^2(Q_T)}^2
   =\int_{\{u_n<0\}}\abssq{\partial_xu_n}
   \leq\frac{C_\delta+\omega(n)}{\delta^2\,n^{1-\delta}}.
\end{split}
\end{equation}
Indeed, $0<\tilde{u}_n\leq1/n$ on $\{u_n<0\}$, so the second inequality
follows by inserting $\tilde{u}_n^{1-\delta}\tilde{u}_n^{\delta-1}$ and applying
\eqref{eq:delta_apriori}. Also, for every $p\in[\delta,3]$,
\[
 \sup_{0\leq t\leq T}\norm{(u_n)_-(t,\cdot)}_{L^p(\Omega;d\nu_\lambda)}^p
 \leq n^{-p}\nu_\lambda(\Omega_n)\longrightarrow0.
\]
In particular $(u_n)_-\to0$ in $L^2(0,T;H^1(J))$ for every bounded
interval $J\Subset\Omega$. These estimates hold on the full time
interval $(0,T)$ and will be used with the energy and
truncation statements below.

\subsection{Compactness}\label{subsection:compactness}
Next, we obtain precompactness of the family $\{u_n\}_{n\in \bbN}$ in the strong topology on 
$L^{p}_{loc}(Q_T;dt \times d\nu_\lambda)$ in Proposition~\ref{prop:evolution_compactness}. We then extract a candidate solution via subsequential limit in 
Proposition~\ref{prop:convergence_from_compactness} and 
record convergences to it.
\begin{proposition}\label{prop:evolution_compactness}
The family of regularized solutions $\{u_n\}_{n \in \bbN}$ is precompact in 
$L^q_{loc}(\overline{Q}_T; dt \times d\nu_\lambda)$ for all $q \in (\delta, 3)$, with $\delta < 1$ from Theorem~\ref{thm:existence}.
\end{proposition} 
The following lemmas serve the proof of Proposition~\ref{prop:evolution_compactness}.
\begin{lemma}[Aubin--Lions]\label{lem:aubin_lions}
The family $\{u_n\}_{n \in \mathbb{N}}$ is precompact in
\[
  L^1\!\big([0, T]; L^r_{loc}([0, \infty); d\nu_\lambda)\big), 
  \qquad r := \begin{cases} 4/3 & \lambda = 0, \\ 3 & \lambda > 0. \end{cases}
\]
\end{lemma}
First, we establish the Banach space framework for Aubin--Lions in our weighted setting. 
\begin{lemma}[Aubin--Lions Triple]\label{lem:abl_triple}
Fix $0 < b < \infty$. With $V^q_L(0, b)$ as defined Subsection~\ref{subsection:notation} and $r$ as 
in Lemma~\ref{lem:aubin_lions}, define
\[
  A(0, b) := \begin{cases} V^3_L(0, b; d\nu_\lambda) \cap V^1_L(0, b; d\nu_\lambda) & \lambda = 0, \\ 
                          V^3_L(0, b; d\nu_\lambda) & \lambda > 0, \end{cases}
\] 
                      
\[ 
 B(0, b) := L^r(0, b; d\nu_\lambda), \qquad C(0, b) := H^{-1}(0, b).
\]
Then $A \cembed B$ compactly and $B \hookrightarrow C$ continuously.
\end{lemma}
\begin{proof}[Proof of Lemma~\ref{lem:abl_triple}]
The Hardy inequalities of Appendix~\ref{appendix:hardy_sobolev} give the compact embedding 
$V^3_L(0, b; d\nu_\lambda) \cembed L^3(0, b; d\nu_\lambda)$ for all $\lambda \geq 0$. 
When $\lambda > 0$, this immediately gives $A \cembed B$, whereas for $\lambda = 0$, additionally apply 
$L^q$ interpolation for $q \in (1, 3)$,
\[
  \norm{f}_{L^q(0, b; d\nu_\lambda)}^q 
  \leq C(q) \big(\norm{f}_{L^1(0, b; d\nu_\lambda)}^{3-q} 
  \norm{f}_{L^3(0, b; d\nu_\lambda)}^{3(q-1)}\big)^{1/2},
\]
Specializing to $q = 4/3$ to yield $A \cembed B$ for $\lambda = 0$. To see $B \hookrightarrow C$, define the inclusion $\iota: B \to H^{-1}(0, b)$ via
$\iota(f)(\phi) := \ip{f}{\phi}_{L^2(d\nu_\lambda)}$. Injectivity of $\iota$ is 
immediate from $L^p$ duality, and by the 1D Sobolev embedding we have
$\abs{\phi(x)} \leq \norm{\phi}_{H^1} x^{1/2}$ for $\phi \in H^1_0(0, b)$ and 
$x \leq b$. Combine with H\"older's inequality to find,
\[ 
  \abs{\iota(f)(\phi)} 
  \leq \begin{cases} 
    \norm{\phi}_{H^1_0(0, b)} \, b^{1/4} \norm{f}_{L^{4/3}(d\nu_0)} 
      & \lambda = 0, \\
    \nu_\lambda(0, b)^{2/3} C_S(0, b) \norm{\phi}_{H^1_0(0, b)} \norm{f}_{L^3(d\nu_\lambda)} 
      & \lambda > 0
  \end{cases}
\]
using conjugate exponents $(4, 4/3)$ and $(3,3/2)$ respectively. These give the continuous inclusion $B \hookrightarrow C$.
\end{proof}
\begin{proof}[Proof of Lemma~\ref{lem:aubin_lions}]
Using Lemma~\ref{lem:abl_triple}, we show that $\{u_n\}_{n \in \bbN}$ is precompact in 
$L^1(0, T; B(0, b))$ for any $b < \infty$. The Aubin--Lions criterion 
(see~\cite[Theorem II.5.16]{boyer}, \cite[Theorem 5]{aubin_lions_by_simon}) requires that we verify
\begin{equation}
  \sup_n \norm{u_n}_{L^1(0, T; A(0, b))} < \infty, 
  \qquad
  \sup_n \norm{\diff{t} u_n}_{L^1(0, T; C(0, b))} < \infty.
\end{equation}
Fix $b < \infty$ and take $n$ large enough that $b_n \geq b$ in~\eqref{eq:omega_n_requirement}. The first criterion 
follows from Proposition~\ref{prop:uniform_estimates} via $L^p(0, T; A) \hookrightarrow L^1(0, T; A)$ for finite $T$.
For the second criterion, keep the boundary term created by extending
$u_n$ by zero. For $\phi\in H^1_0(0,b)$ and $b_n>b$, equation~\eqref{eq:n_problem}
gives
\[
\begin{split}
 \iota(\partial_tu_n)(\phi)
 &=-\int_{a_n}^b \tilde{u}_n\partial_xu_n\partial_x\phi
   -\int_{a_n}^b\abssq{\partial_xu_n}\phi
   +\int_{a_n}^b g_n\tilde{u}_n\phi\\
 &\quad-\frac1n\partial_xu_n(t,a_n)\phi(a_n).
\end{split}
\]
The Sobolev embedding $H^1_0(0,b)\hookrightarrow C^0([0,b])$ and the
continuous inclusion $B(0,b)\hookrightarrow H^{-1}(0,b)$ proved above
therefore yield
\[
\begin{split}
 \abs{\iota(\partial_tu_n)(\phi)}
 &\leq C_b\norm{\phi}_{H^1_0(0,b)}\Bigl(
      \norm{\tilde{u}_n\partial_xu_n}_{L^2(a_n,b)}
      +\norm{\partial_xu_n}_{L^2(a_n,b)}^2\\
 &\qquad+\linfty{g\rho_\lambda}
       \norm{\tilde{u}_n\mathbf1_{(a_n,b)}}_{B(0,b)}
      +\frac1n\abs{\partial_xu_n(t,a_n)}\Bigr).
\end{split}
\]
Here the indicator means extension by zero to $(0,b)$. The reaction
term is bounded using $\abs{\rho_\lambda g_n}\leq\linfty{g\rho_\lambda}$;
this avoids requiring $\nu_0(0,b)$ to be finite.
By~\eqref{eq:shifted_power_bound} with $p=1,3$ and interpolation, the
first three terms are uniformly integrable in time. The last term is
controlled by~\eqref{eq:n_boundary_flux_l1}. Consequently,
\[
 \sup_n\norm{\iota(\partial_tu_n)}_{L^1(0,T;H^{-1}(0,b))}<\infty,
\]
which is the required time-derivative bound in $C(0,b)$.
To extend precompactness to $[0,\infty)$, simply take right endpoints $b_k \to \infty$ and use a standard diagonal argument.
\end{proof}
\begin{proof}[Proof of Proposition~\ref{prop:evolution_compactness}]
For $\lambda > 0$, $\nu_\lambda((0,\infty)) < \infty$ implies by H\"older's inequality,
\[
\nu_\lambda(0, b)^{-2/3} \norm{u_n - u_{n'}}_{L^1((0,T)\times(0,b); dt \times d\nu_\lambda)}
\leq \norm{u_n - u_{n'}}_{L^1(0, T; L^3(0, b; d\nu_\lambda))}.
\]
It follows that precompactness in $L^1(0, T; L^3_{loc}(0, \infty; d\nu_\lambda))$ from 
Lemma~\ref{lem:aubin_lions} transfers directly to $L^1_{loc}(\overline{Q}_T; dt \times d\nu_\lambda)$. For $\lambda = 0$, the 
$L^{p}-L^{q}$ interpolation pair $(\delta, 4/3)$ gives, 
\[
  \norm{f}_{L^1(0, b; d\nu_\lambda)} 
  \leq C(\delta, 4/3) \norm{f}_{L^\delta(0, b; d\nu_\lambda)}^\alpha 
  \norm{f}_{L^{4/3}(0, b; d\nu_\lambda)}^\beta,
\]
with $\alpha = \delta/(4 - 3\delta)$ and $\beta = (1-\delta)/(1 - \delta + \delta/4)$, 
satisfying $\alpha,\beta \in (0, 1)$. The uniform $L^\delta$ bound on $u_n$ follows from~\eqref{eq:delta_apriori},
$\abs{u_n}^\delta\leq \tilde{u}_n^\delta+n^{-\delta}$ on $\Omega_n$, and
\eqref{eq:omega_n_requirement}. Applying the interpolation inequality
to differences and integrating in time, the $L^{4/3}$ precompactness
from Lemma~\ref{lem:aubin_lions} yields precompactness in
$L^1_{loc}(\overline{Q}_T;dt\times d\nu_\lambda)$.
For all $\lambda\geq0$, interpolation using~\eqref{eq:lp_apriori}
and the uniform $L^\delta$ bound extends this to $q\in(\delta,3)$.
\end{proof}
We conclude this section by identifying a candidate weak solution $u$ via subsequential limit. 
\begin{proposition}[Convergence to a Limit Function]\label{prop:convergence_from_compactness}
There exists a nonnegative function $u \colon Q_T \to \bbR$ and a subsequence (not relabeled) of $\{u_n\}_{n \in \mathbb{N}}$ such that:
\begin{enumerate}[label=(\arabic*)]
  \item $u_n \to u$ strongly in $L^{q}_{loc}(Q_T;dt \times d\nu_{\lambda}), \forall q \in (\delta,3)$
  hence $u_n(t,x) \to u(t,x)$ for $(t,x) \in Q_{T}$ a.e.
  \item $u \in L^{\infty}(0,T;L^{q}(\Omega;d\nu_{\lambda}))$ for $q \in [1,3]$, and additionally $q \in [\delta,3]$ when $\lambda = 0$
  ($\delta$ from Theorem~\ref{thm:existence}). 
  \item $\partial_{x}u_n \wto \partial_{x}u$ and $\partial_{x}u_n^2 \wto \partial_{x}u^2$  in $L^2(Q_T)$;
\end{enumerate} 
\end{proposition}
\begin{proof}
Statement (1) is the direct conclusion of Proposition~\ref{prop:evolution_compactness}; the limit is nonnegative by~\eqref{eq:n_maximum_principle}. For (2), Fubini-Tonelli applied to (1) gives a Lebesgue-null $N \subset (0,T)$ 
such that $u_n(t,\cdot) \to u(t,\cdot)$, $\nu_\lambda$-a.e. (same as Lebesgue a.e.) for $t \in (0,T) \setminus N$. 
Using the a.e. convergence of (1), Fatou's lemma at each $t\in (0,T)\setminus N$ with the uniform-in-t bounds from~\eqref{eq:lp_apriori} and~\eqref{eq:delta_apriori} yields 
$\norm{u(t,\cdot)}_{L^q(\Omega; d\nu_\lambda)}^q \leq C_q$ for $q \in \{1,2,3\}$ 
(plus $q = \delta$ when $\lambda = 0$). As $N$ is Lebesgue null in $[0,T]$, this proves (2) for $q \in \{\delta,1,2,3\}$, and 
$L^{p}$ interpolation gives the stated range $q \in [\delta,3]$. For (3), the $m=0,2$ cases of~\eqref{eq:lp_apriori} bound
$\partial_xu_n$ and $\partial_xu_n^2$ in $L^2(Q_T)$, since
$\abssq{\partial_x\abs{u_n}}=\abssq{\partial_xu_n}$ a.e.
Banach--Alaoglu and the strong local convergence in (1) identify their
weak limits as $\partial_xu$ and $\partial_xu^2$ by testing against
smooth compactly supported functions. We also have
\[
 \partial_x(\tilde{u}_n^2)-\partial_x(u_n^2)
     =\frac2n\partial_xu_n\longrightarrow0
       \quad\text{in }L^2(Q_T),
\]
where $\tilde{u}_n=u_n+1/n$ after the zero extension of $u_n$.
\end{proof}

\section{Existence of Weak Solutions}\label{section:existence}
In this section we prove Theorem~\ref{thm:existence} by passing to the limit in the weak form of~\eqref{eq:n_problem}, 
which reads, for $\phi \in C^\infty_c([0, T) \times \Omega)$ and $n$ sufficiently large,
\begin{equation*}
\begin{split}
& -\fintegral{Q_T}u_n \diff{t}\phi\,dt\,d\nu_{\lambda}
  + \fintegral{Q_T} \partial_{x}{\frac{\left(u_n+\frac{1}{n}\right)^2}{2}}\,\partial_{x}{\phi}\,dt\,dx
  + \fintegral{Q_T} \phi\abssq{\partial_{x}{u_n}}\,dt\,dx \\
& = \fintegral{\Omega} \phi(0,x)\,u_{0,n}(x)\,d\nu_{\lambda}
  + \fintegral{Q_T} \phi\,\bigl(u_n+\tfrac{1}{n}\bigr)\,g_n\,dt\,dx.
\end{split}
\end{equation*}
The convergences established in Proposition~\ref{prop:convergence_from_compactness} suffice to pass to the limit for all terms except the quadratic-gradient term. Subsection~\ref{subsection:gradient_convergence} 
treats this term, establishing strong convergence of $\partial_x u_n$ in $L^2_{loc}(Q_T)$ in Corollary~\ref{cor:gradients_converge_locally}; Subsection~\ref{subsection:existence} concludes the proof
of existence.

\subsection{Strong Convergence of Gradients}\label{subsection:gradient_convergence}
Here we obtain strong $L^2_{loc}$ convergence of gradients $\partial_x u_n \to \partial_x u$, up to subsequence. Throughout, $\tilde{u}_n=u_n+1/n>0$ by~\eqref{eq:shifted_positivity}, with the negative part of $u_n$ controlled by~\eqref{eq:negative_part_estimates}. We first record a lemma used in the forthcoming arguments.
\begin{lemma}[Energy Inequalities]
\label{lem:dallaglio_energy_inequalities}
Let $\psi \in C([0, \infty))$ be Lipschitz, nondecreasing, with $\psi(0) = 0$, and 
$\phi(s) = \lintegral{0}{s} z \psi(z)dz$, extending $\psi$ and $\phi$ by zero on $(-\infty,0)$. For any $n \in \mathbb{N}$:
\begin{equation}\label{eq:psi_energy_inequality}
\begin{split}
  &\sup_{t \in [0, T]}\left\{\fintegral{\Omega_n} \frac{\phi(u_n(t, x))}{\rho_\lambda(x)} \, dx 
  + \fintegral{Q_t} \tilde{u}_n u_n \psi'(u_n) \abssq{\partial_{x}u_n}\right\} \\
  &\quad \leq \fintegral{\Omega_n} \frac{\phi(u_{0, n}(x))}{\rho_\lambda(x)} \, dx 
  + \fintegral{Q_T} \abs{g_n} u_n^2 \psi(u_n) + \omega(n).
\end{split}
\end{equation}
For any $v \in C^\infty_C([0, T] \times \Omega_n)$, with $\sigma_v(u) := 
\mathbf{1}_{v > 0} + \frac{1}{u + \tfrac{1}{n}} \mathbf{1}_{v < 0}$,
\begin{equation}\label{eq:signed_energy_inequality}
\begin{split}
  &\fintegral{Q_T} \diff{t} u_n \cdot \sigma_v(u_n)v \frac{dx}{\rho_\lambda(x)} \, dt
  + \fintegral{Q_T} \sigma_v(u_n) (u_n + \tfrac{1}{n}) \partial_{x}u_n \cdot \partial_{x}v \\
  &\quad \leq \fintegral{Q_T} \sigma_v(u_n) (u_n + \tfrac{1}{n}) g_n v.
\end{split}
\end{equation}
\end{lemma}
\begin{proof}
For~\eqref{eq:psi_energy_inequality}, test~\eqref{eq:energy_form_n} with $\phi'(u_n) = u_n \psi(u_n)$ over $\Omega$,
integrate in time from $0$ to $t \leq T$, and use the zero extensions of $\psi,\phi$ and that $\psi(u_n) \geq 0$ along with
the definition of $\omega(n)$~\eqref{eq:omega_define} to see $\fintegral{Q_{T}} u_n (u_n + \tfrac{1}{n}) \psi(u_n) \abs{g_n} 
= \fintegral{Q_T} u_n^2 \psi(u_n) \abs{g_n} + \omega(n)$, yielding the bound. 
For~\eqref{eq:signed_energy_inequality}, test with $\sigma_v(u_n)v=v_+-v_-/\tilde{u}_n$ and integrate 
over $Q_T$. The quadratic-gradient terms cancel on $\{v<0\}$, while on $\{v>0\}$ the remaining term is $v_+\abssq{\partial_xu_n}\geq0$, giving the bound. The identities extend by approximation to the compactly supported Sobolev test functions used below.
\end{proof}
We next show that where the family $\{u_n\}_{n \in \bbN}$ is uniformly small or large pointwise, 
the gradients $\partial_{x}u_n$ are uniformly small in $L^2(Q_T)$. Estimate~\eqref{eq:gradient_sublevel} 
adapts~\cite[Proposition 4.2]{dallaglio_primary} while~\eqref{eq:gradient_suplevel} is a new direct argument.
\begin{proposition}[Gradient Tail Bounds] \label{prop:gradient_tail_bounds}
For any $k \geq 1$ and $\delta \in (0, 1)$ from Theorem~\ref{thm:existence}, $C_\delta$ from~\eqref{eq:delta_apriori},
\begin{equation}\label{eq:gradient_sublevel}
  \fintegral{Q_T \cap \{u_n \leq 1/k\}} \abssq{\partial_{x}u_n} 
  \leq \big(\tfrac{1}{k} + \tfrac{1}{n}\big)^{1-\delta} \frac{C_\delta+\omega(n)}{\delta^2},
\end{equation}
and
\begin{equation}\label{eq:gradient_suplevel}
  \fintegral{\{u_n \geq k\}} \abssq{\partial_{x}u_n} \leq \omega(k) + \omega(n).
\end{equation}
\end{proposition}
\begin{proof}
For~\eqref{eq:gradient_sublevel}, let $\beta = 1-\delta$ and
introduce $(u_n + 1/n)^\beta$ into the integral
$\fintegral{Q_T \cap \{u_{n}\leq \frac{1}{k}\}} \abssq{\partial_{x}u_n}$ to find

\begin{equation*}
\begin{split}
  \fintegral{Q_T \cap \{u_n \leq 1/k\}} \abssq{\partial_{x}u_n}
  &= \fintegral{Q_T \cap \{u_n \leq 1/k\}} \frac{\left(u_{n}+\frac{1}{n}\right)^{\beta}}{\left(u_{n}+\frac{1}{n}\right)^{\beta}}\abssq{\partial_{x}u_n}\\
  &\leq \big(\tfrac{1}{k} + \tfrac{1}{n}\big)^\beta 
  \fintegral{Q_T} \tfrac{\abssq{\partial_{x}(u_n + \tfrac{1}{n})}}{(u_n + \tfrac{1}{n})^\beta}.
\end{split}
\end{equation*}

Now use the identity $(u_n + \tfrac{1}{n})^{-\beta} \abssq{\partial_{x}(u_n + \tfrac{1}{n})} 
= (1 - \beta/2)^{-2} \abssq{\partial_{x}(u_n + \tfrac{1}{n})^{1 - \beta/2}}$, 
and recognize this last integral is already controlled by the estimate~\eqref{eq:delta_apriori}.
Inserting the bound from~\eqref{eq:delta_apriori} here gives the result. For~\eqref{eq:gradient_suplevel}, we use Lemma~\ref{lem:dallaglio_energy_inequalities}
with the choice $\psi(s) = (\tfrac{1}{k} - \tfrac{1}{s}) \mathbf{1}_{s \geq k}$, 
$\phi(s) = (\tfrac{s^2}{2k} - \tfrac{k}{2} - (s - k)) \mathbf{1}_{s \geq k}$, and 
$\psi'(s) = s^{-2} \mathbf{1}_{s \geq k}$. Inserting these choices into
\eqref{eq:psi_energy_inequality} and dropping nonnegative terms gives,
\[
\begin{aligned}
  \fintegral{Q_T\cap\{u_n \geq k\}} \abssq{\partial_{x}u_n}
  &\leq \fintegral{\{u_{0, n} \geq k\}} \frac{\phi(u_{0, n}(x))}{\rho_\lambda(x)} \, dx \\
  &\quad + \fintegral{\{u_n \geq k\}} \abs{g_n(x)} u_n \big(\tfrac{u_n}{k} - 1\big) + \omega(n).
\end{aligned}
\]
The first term above is $\omega(k)$, uniformly in $n$, since
$0\leq\phi(s)\leq s^2/(2k)$ and $u_{0,n}\to u_0$ in
$L^2(\Omega;d\nu_\lambda)$ by~\eqref{eq:n_problem_data_requirement}.
The second term is $\omega(k)+\omega(n)$ by the uniform estimates~\eqref{eq:lp_apriori}
and the assumptions on $g_n$:
\[
 \fintegral{Q_T\cap\{u_n\geq k\}}\abs{g_n}u_n
       \bigl(\tfrac{u_n}{k}-1\bigr)
 \leq \frac{\linfty{g\rho_\lambda}}{k}
       \normsq{u_n}_{L^2(Q_T;dt\times d\nu_\lambda)}
 =\omega(k)+\omega(n),
\]
completing the proof. 
\end{proof}

By Proposition~\ref{prop:gradient_tail_bounds}, the gradients $\partial_{x}u_n$ are 
vanishing in $L^2$, jointly in $k$ and $n$, when $u_n$ lie outside the range $[1/k, k]$. Accordingly, introduce the double truncation
\begin{equation}\label{defn:truncation_define}
  T_k(s) := \begin{cases} 
    1/k & s < 1/k, \\
    s & 1/k \leq s \leq k, \\
    k & s > k,
  \end{cases}
\end{equation}
with complement
\begin{equation}\label{eq:identity_cutoff_decomp}
  G_k(s) := s - T_k(s) = (s - k)_+ - (1/k - s)_+.
\end{equation}
We now prove
local, strong convergence $\partial_{x}T_k u_n \to \partial_{x}T_k u$ in $L^2(Q_T)$ using the steps of~\cite[Proposition 6.2]{dallaglio_primary}. In the sequel, we will often write $T_k u_n$ in lieu of $T_k(u_n)$. 
\begin{proposition}[Gradients of Truncations Converge Strongly]\label{prop:trunc_gradients_converge}
For any fixed $k > 1$ and nonnegative $\zeta \in C^\infty_C([0, T) \times \Omega)$,
\begin{equation}\label{eq:trunc_gradients_converge}
  \fintegral{Q_T} \zeta \abssq{\partial_{x}\btk u_n - \partial_{x}\btk u} \longrightarrow 0 
  \quad \text{as } n \to \infty.
\end{equation}
\end{proposition}
\begin{proof}
Let $\phi_+(s) := \exp(\alpha s_+) - 1$ for any $\alpha \geq 1$ (further 
constraints on $\alpha$ appear in the sequel). Introduce the time regularization, for any $\eta \in \bbN$,
\begin{equation}\label{eq:u_eta_define}
  (\btk u)_\eta(t, x) := e^{-\eta t} \btk u_{0, j(\eta)}(x) 
  + \eta \lintegral{0}{t} e^{-\eta(t - s)} \btk u(s, x) \, ds,
\end{equation}
which solves $\tfrac{1}{\eta} \diff{t} (\btk u)_\eta + (\btk u)_\eta = \btk u$ 
with initial data $(\btk u)_\eta(0, x) = \btk u_{0,j(\eta)}$ (cf.~\cite{landes_first,landes_second}). Here $j(\eta)\to\infty$ is chosen sufficiently slowly that
$\eta^{-1}\normsq{\partial_xu_{0,j(\eta)}}_{L^2(\Omega)}\to0$; this ensures strong $L^2(Q_T)$ convergence of the regularized spatial gradients without an initial-gradient assumption.
Fix nonnegative $\zeta \in C^\infty_C((0, \infty))$ 
and form the test function
\[
  v := \zeta \, \phi_+(\btk u_n - (\btk u)_\eta).
\]
With this choice, the truncated part of the time term is handled by the
ordinary primitive $\Psi_+(s)=\int_0^s\phi_+(r)\,dr$.
Insert this to~\eqref{eq:signed_energy_inequality} of Lemma~\ref{lem:dallaglio_energy_inequalities}, giving:
\begin{equation}\label{eq:three_terms_start}
  \underbrace{\fintegral{Q_T} (\diff{t} u_n) v\,dt\,d\nu_\lambda}_{\mathcal{T}_1} 
  + \underbrace{\fintegral{Q_T} (u_n + \tfrac{1}{n}) \partial_{x}u_n \cdot \partial_{x}v}_{\mathcal{T}_3} 
  \leq \underbrace{\fintegral{Q_T} g_n (u_n + \tfrac{1}{n}) v}_{\mathcal{T}_2}.
\end{equation}
Since $(T_ku)_\eta\geq1/k$, both $v$ and $\partial_xv$ vanish a.e. on $\{u_n\leq1/k\}$. Thus the negative set contributes nothing to these three terms, and the integrals may remain over $Q_T$.
To prove the Proposition, we now prove 3 Claims, and synthesize the outcome: ultimately we extract gradient information from $\mathcal{T}_3$ and prove other terms have 
favorable sign or vanish. For two-index sequences $a_{\eta, n}$, we write
$a_{\eta, n} \equiv \omega^\eta(n)$ to mean $\limsup_{n \to \infty} |a_{\eta, n}| = 0$ 
for each fixed $\eta$. The combined errors below are used in the successive-limit sense $\lim_{\eta\to\infty}\limsup_{n\to\infty}|\mathrm{error}_{\eta,n}|=0$. Following~\cite{dallaglio_primary}, we verify three claims:
\begin{description}
  \item[Claim 1.] $\mathcal{T}_1 \geq \omega(\eta) + \omega^\eta(n)$,\ i.e. $\mathcal{T}_1$ is asymptotically nonnegative;
  \item[Claim 2.] $\mathcal{T}_2 = \omega(n) + \omega^\eta(n) + \omega(\eta)$,\ i.e. $\mathcal{T}_2$ is asymptotically zero;
  \item[Claim 3.] $\mathcal{T}_3 \geq \fintegral{Q_T\cap\{1/k\leq u_n\leq k\}} \zeta \abssq{\partial_x(\btk u_n-(\btk u)_\eta)_+} 
       + \omega(\eta) + \omega^\eta(n)$,\ i.e. $\mathcal{T}_3$ asymptotically dominates squared gradient distance.
\end{description}
The proofs of Claims 1-3 are presented in Appendix~\ref{appendix:dall_aglio}: proofs of Claims 1 and 2 closely follow~\cite[Proposition 6.2]{dallaglio_primary}, modified to accommodate $\rho_{\lambda}$,
while verification of Claim 3 keeps the shifted coefficient $\tilde{u}_n>0$ throughout.
When input to estimate~\eqref{eq:three_terms_start}, the conclusions of Claims 1-3 yield 
\begin{equation}
  \fintegral{\{1/k \leq u_n \leq k\}} \zeta \abssq{\partial_{x}(\btk u_n - (\btk u)_\eta)_+}
  \leq \omega(n) + \omega^\eta(n) + \omega(\eta).
\end{equation}
For the negative difference, retain $\phi_-(s) := 1 - \exp(\alpha s_-)$ and
the test $v=\zeta\phi_-(\btk u_n-(\btk u)_\eta)/\btk u_n$
(cf.~\cite{dallaglio_secondary}). Using the negative branch of~\eqref{eq:signed_energy_inequality} yields the analogous bound on
$\abssq{\partial_{x}(\btk u_n - (\btk u)_\eta)_-}$.
The gradients of the positive and negative parts have disjoint supports a.e., producing the identity
\[ \partial_x(\btk u_n-(\btk u)_\eta)
 =\partial_x(\btk u_n-(\btk u)_\eta)_+
  -\partial_x(\btk u_n-(\btk u)_\eta)_-. \]
Outside $\{1/k\leq u_n\leq k\}$, the remaining squared gradient is $\abs{\partial_x(T_ku)_\eta}^2$; its integral tends to zero, first in $n$ and then in $\eta$, by strong convergence $\partial_x(T_ku)_\eta\to\partial_xT_ku$ and the fact that $\partial_xT_ku=0$ a.e. on $\{u\notin(1/k,k)\}$.
Combining the bounds on the last two terms gives
\[ 
\begin{split} 
  & \fintegral{Q_T}  \abssq{\partial_{x}{}(\btk u_n - (\btk u)_\eta)}\zeta = \fintegral{Q_T} (\abssq{\partial_{x}{}(\btk u_n - (\btk u)_\eta)_{+}} \\ 
  & + \abssq{\partial_{x}{}(\btk u_n-(\btk u)_\eta)_{-}})\zeta \leq \omega(\eta) + \omega^{\eta}(n).
\end{split}
 \]
Since we also have $\partial_{x}(\btk u)_\eta \to \partial_{x}\btk u$ strongly in $L^2(Q_T)$ as 
$\eta \to \infty$ by construction~\cite{landes_second}, we pass (by triangle inequality)
from $u_{\eta}$ to $u$ in the last bound, implying the main result~\eqref{eq:trunc_gradients_converge} for $\zeta \in C^\infty_C((0, \infty))$. 
For $\zeta\in C^\infty_C([0,T)\times\Omega)$, dominate $\zeta$ by a constant multiple of a nonnegative spatial cutoff. The preceding estimates already integrate from $0$ to $T$, completing the proof.
\end{proof}
These next corollaries facilitate the proofs of both Theorems~\ref{thm:existence} and~\ref{thm:regularity}.
\begin{corollary}[Gradients Converge a.e.]\label{cor:gradients_almost_everywhere}
Up to a subsequence,
\begin{equation}
  \partial_{x}u_n \to \partial_{x}u \quad \text{a.e.}
\end{equation}
\end{corollary}
\begin{proof}
For any $k \geq 1$, decompose
\begin{equation}
  \partial_xu_n-\partial_xu
  =\partial_x(T_ku_n-T_ku)+\partial_xG_k(u_n)-\partial_xG_k(u).
\end{equation}
For each fixed $k$, the first term converges to zero in $L^2((0,S)\times J)$ for $S<T$ and $J\Subset\Omega$, by Proposition~\ref{prop:trunc_gradients_converge}.
For the second term, the gradient tail bounds~\eqref{eq:gradient_sublevel}--\eqref{eq:gradient_suplevel} give
$\lim_{k\to\infty}\limsup_{n\to\infty}\norm{\partial_xG_k(u_n)}_{L^2(Q_T)}=0$.
Finally, $\partial_xG_k(u)\to0$ in $L^2(Q_T)$ by dominated convergence, since $u\geq0$ and $\partial_xu=0$ a.e. on $\{u=0\}$.
Chebyshev's inequality therefore gives local convergence in measure of the full gradients. A diagonal subsequence on an exhaustion of $Q_T$ gives the stated a.e. convergence.
\end{proof}
\begin{corollary}[Strong, Local Gradient Convergence]\label{cor:gradients_converge_locally}
Up to a subsequence,
\begin{equation}\label{eq:gradients_converge_locally}
  \partial_{x}u_n \to \partial_{x}u \quad \text{in } L^2_{loc}(Q_T).
\end{equation}
The convergence also holds in $L^2((0,S)\times J)$ for every $S<T$ and $J\Subset\Omega$ (i.e. including near $t =0$).
\end{corollary}
\begin{proof}
By Vitali's convergence theorem, it suffices to show that $\abssq{\partial_{x}u_n}$ is 
locally uniformly integrable on $Q_T$, since a.e.\ convergence is established in 
Corollary~\ref{cor:gradients_almost_everywhere}. Fix $S<T$, $J\Subset\Omega$. For measurable $A\subset(0,S)\times J$ and $k>1$,
\[
  \fintegral{A} \abssq{\partial_{x}u_n} 
  = \fintegral{A} \abssq{\partial_{x}\btk u_n} 
  + \fintegral{A \cap \{u_n \in (-\infty,\frac{1}{k}) \cup (k,\infty)\}} \abssq{\partial_{x}u_n}.
\]
The second term is $\omega(k)+\omega(n)$ via Proposition~\ref{prop:gradient_tail_bounds}, so given $\epsilon>0$ we first fix $k$ sufficiently large and then $n$ sufficiently large to make it less than $\epsilon/3$.
Meanwhile, Proposition~\ref{prop:trunc_gradients_converge} gives for all $k$,
\[
  \fintegral{A} \abssq{\partial_{x}\btk u_n} = \fintegral{A} \abssq{\partial_{x}\btk u} + \omega(n) 
  \leq \fintegral{A} \abssq{\partial_{x}u} + \omega(n),
\]
with the error uniform over $A$. Absolute continuity of the Lebesgue integral yields $\eta>0$ such that $|A|<\eta$ makes the limiting integral less than $\epsilon/3$; the remaining error is less than $\epsilon/3$ for $n$ large. The finitely many remaining indices are handled by absolute continuity of their integrals. Uniform integrability follows, establishing~\eqref{eq:gradients_converge_locally}, including near $t=0$.
\end{proof}
\subsection{Existence}\label{subsection:existence}
Now we are in position to conclude the proof of Theorem~\ref{thm:existence} by passing to the limit in the weak form of~\eqref{eq:n_problem}.

\begin{proof}[Proof of Theorem~\ref{thm:existence}]
Take the subsequential limit $u$ of Proposition~\ref{prop:convergence_from_compactness}. That $u$ satisfies the regularity claims
\[
  L^{\infty}(0,T;L^q(\Omega;d\nu_{\lambda})) \text{ for } q \in [\delta,3]
  \quad\text{and}\quad \partial_{x}u \in L^2(Q_T)
\]
are both concluded in Proposition~\ref{prop:convergence_from_compactness}.
To see that $u^2 \in L^{3/2}(0,T; V_0^{3/2}(\Omega))$, note each $u_n^2(t,\cdot)$ extended by zero to 
$\Omega$  lies in $V_0^{3/2}(\Omega)$ since $u_n|_{\partial\Omega_n} = 0$ due to Proposition~\ref{prop:n_are_classical}, and 
$\{u_n^2\}_{n \in \bbN}$ is uniformly bounded there by estimates~\eqref{eq:lp_apriori} ($m=2$ case).
Hence, Banach--Alaoglu for $\{u_n^2\}_{n \in \bbN}$, uniformly bounded in the reflexive space
$L^{3/2}(0,T; V_0^{3/2}(\Omega))$,
plus the a.e.\ convergence of Proposition~\ref{prop:convergence_from_compactness}, prove
$u^2 \in L^{3/2}(0,T; V_0^{3/2}(\Omega))$. 
To see $u$ satisfies the distributional equation~\eqref{eq:weak_form}, use that $u_n$ solves~\eqref{eq:n_problem} classically 
so that for $\phi \in C^{\infty}_c([0,T) \times \Omega)$ and $n$ sufficiently large,
\begin{equation}\label{eq:weak_form_n}
\begin{split} 
& \underbrace{-\fintegral{Q_T}u_n \diff{t}\phi\,dt\,d\nu_\lambda}_{\mc{I}_1}
 + \underbrace{\fintegral{Q_T} \partial_x\frac{\tilde{u}_n^2}{2}\,\partial_x\phi}_{\mc{I}_2}
 + \underbrace{\fintegral{Q_T}\phi\abssq{\partial_xu_n}}_{\mc{I}_3}\\
&=\underbrace{\fintegral{\Omega}\phi(0,x)u_{0,n}\,d\nu_\lambda}_{\mc{I}_4}
 +\underbrace{\fintegral{Q_T}\phi \tilde{u}_ng_n}_{\mc{I}_5}.
\end{split}
\end{equation}
We are done once each term of~\eqref{eq:weak_form_n} converges to its corresponding term in~\eqref{eq:weak_form}. 
Convergence of $\mathcal{I}_1$ is due to the strong convergence $u_n \to u$ in 
$L^p_{loc}(\overline{Q}_T;dt\times d\nu_\lambda)$ for $p\in(\delta,3)$ from Proposition~\ref{prop:evolution_compactness}, applied on $\support\phi$.
Convergence of $\mathcal{I}_2$ is due to weak convergence $\partial_{x}(u_n + \tfrac{1}{n})^2 \wto \partial_{x}u^2$, also 
from Proposition~\ref{prop:convergence_from_compactness}. For convergence of $\mathcal{I}_3$, 
apply Cauchy-Schwarz and conclude due to Corollary~\ref{cor:gradients_converge_locally} and $L^2(Q_T)$ boundedness of $\partial_{x}{u_n}$ and $\partial_{x}u$ that
\[ \begin{split}
   & \abs{\fintegral{Q_T} \phi(\normsq{\partial_{x}u_n} - \normsq{\partial_{x}u})}\\  
  & \leq \norm{\phi}_{L^\infty(Q_T)}\norm{\mathbf{1}_{\support \phi}(\partial_{x}u_n-\partial_{x}u)}_{L^2(Q_T)}\norm{\partial_{x}u_n+\partial_{x}u}_{L^2(Q_T)}  \xrightarrow{n \to \infty} 0.
  \end{split}
\]
Convergence of $\mathcal{I}_4$ is handled by convergence of initial data, imposed in~\eqref{eq:n_problem_data_requirement}. 
Finally, to see convergence of $\mathcal{I}_5$, observe 
\[
  \abs{\fintegral{Q_T} \phi (u_n + \tfrac{1}{n}) g_n - \fintegral{Q_T} \phi u g}
  \leq I + II + \omega(n),
\]
where
\[
  I := \fintegral{Q_T} \abs{\phi} \abs{u_n - u} \abs{g}, 
  \qquad II := \fintegral{Q_T} \abs{\phi} \abs{u_n} \abs{g_n - g}.
\]
We introduce $\rho_{\lambda}$ and apply H\"older's inequality with exponents $(3,3/2)$ to find 
\[
  I \leq \linfty{\phi} \linfty{\rho_\lambda g} 
  \norm{\mathbf{1}_{\support \phi} (u_n - u)}_{L^1(Q_T; dt \times d\nu_\lambda)} 
  \xrightarrow{n \to \infty} 0,
\]
\[
  II \leq \linfty{\phi} \norm{u_n}_{L^3(Q_T; dt \times d\nu_\lambda)} 
  \norm{\mathbf{1}_{\support \phi} \rho_\lambda (g_n - g)}_{L^{3/2}(Q_T; dt \times d\nu_\lambda)} 
  \xrightarrow{n \to \infty} 0.
\]
Convergence of $I$ follows from $L^{p}_{loc}$ convergence $u_n \to u$ from Proposition~\ref{prop:convergence_from_compactness} 
while the estimate of $II$ follows from $L^{p}_{loc}$ convergence $g_n \rho_\lambda \to g \rho_\lambda$ from~\eqref{eq:n_problem_data_requirement} and uniform bounds on $u_n$ from Proposition~\ref{prop:uniform_estimates}.  
\end{proof}

\section{Smoothing}\label{section:regularity}
This section culminates in a proof of Theorem~\ref{thm:regularity}, the positive-time parabolic H\"older and spatial Lipschitz regularity of the weak solution $u$ found in Theorem~\ref{thm:existence}. 
We first establish the driving estimate in Proposition~\ref{prop:benilan_crandall}, namely, the B\'enilan--Crandall inequality $\partial_t u \geq -u/K_\theta(t)$ in the sense of distributions in $Q_T$.
This inequality and the consequent regularity of $u$ are proved by transferring estimates from the regularized solutions $u_n$ along with standard parabolic interpolation arguments.

\subsection{B\'enilan--Crandall Inequality}\label{subsection:bc}  
\begin{proposition}[B\'enilan--Crandall for $u_n$]\label{prop:bc_regularized}
With $\theta$ and $K_{\theta}$ as in~\eqref{eq:defn_theta} and~\eqref{eq:defn_k_theta}, the classical solution $u_n$ to Problem~\eqref{eq:n_problem}, 
extended by zero to $Q_T$ (cf. Proposition~\ref{prop:n_are_classical}) satisfies classically,
\begin{equation} \label{eq:bc_for_n}
            \diff{t}u_{n}(t,x) \geq -\frac{1}{K_{\theta}(t)}\left(u_{n}+\frac{1}{n}\right),\ \forall\ (x,t) \in Q_{T}.
\end{equation} 
\end{proposition}
\begin{proof}
It suffices to prove~\eqref{eq:bc_for_n} on $Q_{T,n} = (0, T) \times \Omega_n$: outside $Q_{T,n}$, $u_n \equiv 0$ and $\partial_t u_n = 0$ (both classically), including at the boundary of $Q_{T,n}$ (Proposition~\ref{prop:n_are_classical}), making the inequality trivial as the right-hand side is nonpositive.
\textbf{Case $\bm{\theta = 0}$}.\ Write $\tilde{u}_n =u_n + 1/n$ as in~\eqref{eq:defn_tilde_un}, which by direct substitution solves
 $\partial_t \tilde{u}_n = \rho_{\lambda} \tilde{u}_n \partial_x^2 \tilde u_n + \rho_{\lambda} g_n \tilde u_n$ 
classically on $Q_{T,n}$ with boundary data $\tilde u_n(t, x) = 1/n > 0$ for $x \in \partial \Omega_n$. By~\eqref{eq:shifted_positivity}, $\tilde u_n \geq 1/n > 0$ on $\overline{Q_{T,n}}$.
Define the time-scaling of $\tilde{u}_n$ by $K_0(t) = t$, 
\[ z_n(t,x) \defeq t \cdot \diff{t}\tilde{u}_n(t,x) + \tilde{u}_n(t,x)\] 
which, by parabolic regularity at positive times, is a classical solution of the following linear problem and extends continuously to $t=0$:
\begin{equation}
\begin{split}\label{eq:z_benilan_pde}
    z_{n}(0,x) &\defeq u_{0,n}(x)+\frac{1}{n},\ x \in \Omega_n\\
    \diff{t}z_n(t,x) &= \rho_{\lambda}(x)\tilde{u}_n\laplacian{z_n} + (\rho_{\lambda}\laplacian \tilde{u}_n + \rho_{\lambda}g_{n})z_n + \rho_{\lambda}g_{n}\tilde{u}_n,\ (t,x) \in Q_{T,n} \\ 
    z_{n}(t,a_n) &= z_{n}(t,b_n) = \frac{1}{n}.
\end{split}
\end{equation}
The lateral-boundary condition for $z_n$ is valid since $\diff{t}u_n$ and $u_n$ are classical at $\partial Q_{T,n}$. Due to $\tilde{u}_n \geq \frac{1}{n}$, \eqref{eq:z_benilan_pde} is uniformly parabolic, 
thus enjoying comparison principle in the solution $z_n$, since its zeroth-order coefficient is bounded. As $g_n \geq 0$, and $z_n(t,x) \geq 0$ on the parabolic boundary, we find 
$z_{n} \geq 0$ in $Q_{T,n}$, which gives exactly desired result, i.e.~\eqref{eq:bc_for_n} on $Q_{T,n}$, since
\[ z_{n}(t,x) \geq 0 \iff \diff{t}\tilde{u}_n(t,x) \geq -\frac{\tilde{u}_n(t,x)}{t}\ \forall\ (t,x) \in Q_{T,n}. \]
\textbf{Case $\bm{\theta > 0}$}.\ 
Introduce the rescaled time $s(t) =  K_{\theta}(t) = \frac{1}{\theta}\frac{\exp(\theta t)-1}{\exp(\theta t)}$, 
which is strictly increasing from $s(0) = 0$ to $s(\infty) = 1/\theta$ and satisfies $\partial_s = e^{\theta t} \partial_t$. 
Now define $\tilde{u}_{n,\theta} = e^{\theta t} \tilde u_n$ and $g_{n,\theta} 
= g_n + \theta/\rho_{\lambda}$. A direct calculation shows that, in rescaled coordinates $(s,x)$ on $Q_{s(T), n}= (0,s(T)) \times \Omega_n$,
$\tilde u_{n,\theta}$ satisfies
\begin{equation}\label{eq:rescaled_n_theta_problem}
  \partial_s \tilde u_{n,\theta} 
  = \rho_{\lambda} \tilde u_{n,\theta} \partial_x^2 \tilde u_{n,\theta} 
  + \frac{\rho_{\lambda} g_{n,\theta}}{1 - \theta s} \tilde u_{n,\theta},\ \forall\ (s,x) \in Q_{s(T), n}
\end{equation}
with $\tilde u_{n,\theta}(0, x) = u_{0,n}(x)+\frac{1}{n}$ and $\tilde u_{n,\theta}
(s, x) = (1 - \theta s)^{-1}/n$ on $\partial \Omega_n$. The reaction coefficient 
is nonnegative since $\rho_{\lambda} g_{n,\theta} \geq 0$ and 
$1 - \theta s = e^{-\theta t} > 0$. Thus, comparison principle implies that $\tilde{u}_{n,\theta} \geq \frac{1}{n} > 0$ in $Q_{s(T), n}$. 
Now, define $z_{n,\theta}(s, x) := s \cdot \partial_s \tilde u_{n,\theta}(s,x) 
+ \tilde u_{n,\theta}(s,x)$, differentiate in $s$ and use~\eqref{eq:rescaled_n_theta_problem} 
to find that $z_{n,\theta}$ solves the following linear problem,
\begin{equation}
\begin{split}
z_{n,\theta}(0,x) &=  u_{0,n}(x)+\frac{1}{n},\ x \in \Omega_{n}\\
z_{n,\theta}(s,x) &=  \frac{1}{n(1-\theta s)^2},\  s \in (0,s(T)),\ x\in \partial \Omega_{n} \\
\diff{s}z_{n,\theta} &=  \rho_{\lambda}\tilde{u}_{n,\theta}\laplacian z_{n,\theta} + \rho_{\lambda}z_{n,\theta} \laplacian\tilde{u}_{n,\theta} + \rho_{\lambda}\frac{g_{n,\theta}}{1-\theta s}z_{n,\theta}\\
& + \frac{\rho_{\lambda} g_{n,\theta}\tilde u_{n,\theta}}{(1-\theta s)^2},\ \forall\ (s,x) \in Q_{s(T), n}.
\end{split}
\end{equation}
We may compare $z_{n,\theta}$ with zero from below as the equation is uniformly parabolic, its zeroth-order coefficient is bounded, and its source is nonnegative. At the boundary $z_{n,\theta}(0, x) = u_{0,n}+\frac{1}{n} 
\geq 0$ and $z_{n,\theta}(s, a_n) = z_{n,\theta}(s, b_n) > 0$. Consequently, $z_{n,\theta} \geq 0$ on $Q_{s(T), n}$,
and transforming back to the original time,
\[\partial_s \tilde u_{n,\theta} \geq -\frac{\tilde u_{n,\theta}}{s} \iff 
  e^{\theta t}(\partial_t \tilde u_n + \theta \tilde u_n) 
  \geq -\theta\frac{\exp(\theta t)}{\exp(\theta t)-1}\tilde{u}_n,\ \forall\ (t,x) \in Q_{T,n}.
\]
Simplifying the last expression gives exactly~\eqref{eq:bc_for_n}.
\end{proof}  

\begin{proof}[Proof of Proposition~\ref{prop:benilan_crandall}]
Let $\phi \in C^\infty_c(Q_T)$ with $\phi \geq 0$. Multiplying~\eqref{eq:bc_for_n} by $\phi$ and integrating by parts in $t$ gives 
\begin{equation*}
  -\int_{Q_T} u_n \, \partial_t \phi \, dt \, dx 
  \geq -\int_{Q_T} \frac{u_n + \tfrac{1}{n}}{K_\theta(t)} \phi \, dt \, dx.
\end{equation*}
We pass to the limit on both sides as $n \to \infty$ by the strong local convergence $u_n \to u$ in $L^p_{loc}(Q_T; dt \times d\nu_\lambda)$ from Proposition~\ref{prop:convergence_from_compactness} (on $\support \phi$)
with the $\frac{1}{n}$ term vanishing. The result is precisely~\eqref{eq:bc_for_u} in the distributional sense on $Q_T$.
\end{proof}

\subsection{Improved Regularity}\label{subsection:regularity}
Here, we conclude the proof of Theorem~\ref{thm:regularity}. Inequality~\eqref{eq:bc_for_n} yields uniform estimates on the $L^{\infty}$ norms of $u_n$ 
and $\partial_x u_n$ at positive times, which transfer to $u$ in the limit.

\begin{proposition}[Estimates for Positive Times]
\label{prop:improved_estimates}
For $0 < t_0 \leq T$, define
\[
  B_m(t_0,g,u_0) := \frac{(m+1)^2}{4m}\left(\tfrac{1}{K_\theta(t_0)} + \norm{g \rho_{\lambda}}_{L^\infty(\Omega)}\right)C_m(u_0), 
  \qquad m \in \{1, 3\},
\]
with $C_1(u_0)$ and $C_3(u_0)$ from Proposition~\ref{prop:uniform_estimates}. Then, for all $n$ large enough,
\begin{equation}
\begin{split} 
 & \sup_{t_0 \leq t \leq T} \normsq{\partial_{x}u_n(t, \cdot)}_{L^2(\Omega)} \leq B_1(t_0,g,u_0) + \omega(n)\\ 
 & \sup_{t_0 \leq t \leq T} \normsq{\partial_{x}u_n^2(t, \cdot)}_{L^2(\Omega)} \leq B_3(t_0,g,u_0) + \omega(n).
\end{split}
\end{equation}
\end{proposition}
\begin{proof}
Insert the equation for $\partial_t u_n$~\eqref{eq:n_problem} into the B\'enilan--Crandall inequality~\eqref{eq:bc_for_n} and divide through by $\rho_{\lambda}$ to obtain, 
\[
  \left(u_n+\frac{1}{n}\right)\partial_x^2 u_n + g_n \left(u_n+\frac{1}{n}\right) \geq -\frac{1}{\rho_{\lambda}}\left(u_{n}+\frac{1}{n}\right)/ K_\theta(t).
\]
Testing the inequality with 1 and then with $\tilde{u}_n^2/2$, where $\tilde{u}_n=u_n+1/n>0$, respectively yields, after integration by parts, the pointwise bounds for any $t \in [t_0, T]$,
\begin{align*}
 \normsq{\partial_xu_n(t,\cdot)}_{L^2(\Omega)}
 &\leq \left(\frac{1}{K_\theta(t)}+\linfty{g\rho_\lambda}\right)
       \int_{\Omega_n}\tilde{u}_n\,d\nu_\lambda
       +\frac1n[\partial_xu_n]_{a_n}^{b_n},\\
 3\int_{\Omega_n}\tilde{u}_n^2\abssq{\partial_xu_n}
 &\leq \left(\frac{1}{K_\theta(t)}+\linfty{g\rho_\lambda}\right)
       \int_{\Omega_n}\tilde{u}_n^3\,d\nu_\lambda
       +\frac1{n^3}[\partial_xu_n]_{a_n}^{b_n}.
\end{align*}
Here we used $\linfty{g_n\rho_{\lambda}} \leq \linfty{g \rho_{\lambda}}$ from~\eqref{eq:n_problem_data_requirement}. By~\eqref{eq:n_boundary_one_sided}, the boundary terms are bounded above by $2/(n^{p+1}r_n)=\omega(n)$ for $p=1,3$, uniformly in time. Since $u_n^2\leq \tilde{u}_n^2+1/n^2$,
\[
 \normsq{\partial_xu_n^2(t,\cdot)}_{L^2(\Omega)}
 \leq 4\int_{\Omega_n}\tilde{u}_n^2\abssq{\partial_xu_n}
       +\frac4{n^2}\normsq{\partial_xu_n(t,\cdot)}_{L^2(\Omega)}.
\]
Inserting the shifted power bounds~\eqref{eq:shifted_power_bound} from Proposition~\ref{prop:uniform_estimates} and using monotonicity $K_\theta(t) \geq K_{\theta}(t_0)$ for all 
$\theta \geq 0$ gives the result.
\end{proof}
Recall $\mc{D} \Subset Q_{(t_0, T)}$ means that $\mc{D}$ is compactly contained in the cylinder $Q_{(t_0, T)}$, thus having positive distance to the boundary.
\begin{proposition}[Uniform Interior Estimates]\label{prop:interior_bounds_n}
If $\mc{D} \Subset Q_{(t_0, T)}$, there exist constants $\mc{M}_0(\mc{D})$ and $\mc{M}_1(\mc{D})$ such that for all sufficiently large $n$,
\begin{equation}\label{eq:interior_bounds_n}
  \norm{u_n}_{L^\infty(\mathcal{D})} \leq \mc{M}_0(\mc{D}) + \omega(n),\qquad 
  \norm{\partial_x u_n}_{L^\infty(\mathcal{D})} \leq \mc{M}_1(\mc{D}) + \omega(n).
\end{equation}
Let $d:=\inf\{x:(t,x)\in\mc{D}\}>0$, $b:=\sup\{x:(t,x)\in\mc{D}\}<\infty$, and
\[
 N(\mc{D},g):=\frac{1/K_\theta(t_0)+\linfty{g\rho_\lambda}}{\rho_\lambda(d/2)}.
\]
One may take
\[\mc{M}_0(\mc{D}) = \sqrt{(b+d/2)B_1}, \qquad \mc{M}_1(\mc{D}) = \frac{4 \mc{M}_0(\mc{D})}{d} + \frac{N(\mc{D}, g) d}{4} \] 
with $B_1$ from Proposition~\ref{prop:improved_estimates}. 
\end{proposition}
\begin{proof}
Fix $\mc{D}$ and take $n$ large enough that $J:=(d/2,b+d/2)\Subset\Omega_n$. Substituting~\eqref{eq:n_problem} into~\eqref{eq:bc_for_n}, dividing through by $\rho_\lambda(u_n+1/n)>0$ using~\eqref{eq:shifted_positivity}, and using monotonicity of $\rho_\lambda$ gives
\[
  \partial_x^2u_n(t,x)\geq -\frac{1}{K_\theta(t)\rho_\lambda(x)}-g_n(x)
  \geq-N(\mc{D},g),\qquad (t,x)\in[t_0,T]\times J.
\]
For the bound on $\norm{u_n}_{L^{\infty}(\mc{D})}$ in~\eqref{eq:interior_bounds_n}, the zero extension and Cauchy--Schwarz give, for $(t,x)\in[t_0,T]\times J$,
\[
 |u_n(t,x)|=\left|\int_0^x\partial_xu_n(t,y)\,dy\right|
 \leq\sqrt{x}\norm{\partial_xu_n(t,\cdot)}_{L^2(\Omega)}
 \leq\mc{M}_0(\mc{D})+\omega(n).
\]
For the bound on $\norm{\partial_xu_n}_{L^{\infty}(\mc{D})}$, apply Lemma~\ref{lem:lip_from_semiconvex} on $J$ with $M=\mc{M}_0+\omega(n)$ and $\delta\uparrow d/2$, giving the stated formula for $\mc{M}_1$.
\end{proof}
We are now in position to prove Theorem~\ref{thm:regularity}.
\begin{proof}[Proof of Theorem~\ref{thm:regularity}]
It suffices to prove the claims for a fixed cylinder $\mc{D} = (t_1, t_2) \times (a,b) \Subset Q_{(t_0,T)}$, where $0<t_0<t_1<t_2<T$,
and estimate norms on a general compactly contained subset by covering with cylinders. Take $n$ large enough that $\mc{D} \Subset (t_0,T)\times \Omega_n$.
The strong and a.e. convergences $u_n \to u$ of Proposition~\ref{prop:convergence_from_compactness} and $\partial_{x}u_n \to \partial_{x}u$ from Corollary~\ref{cor:gradients_almost_everywhere} and Corollary~\ref{cor:gradients_converge_locally}
imply that the uniform estimates of Proposition~\ref{prop:interior_bounds_n} transfer to $u$ and $\partial_{x}u$ in the limit,
giving both $\norm{u}_{L^\infty(\mc{D})} \leq \mc{M}_0(\mc{D}, g)$ and $\norm{\partial_{x}u}_{L^\infty(\mc{D})} \leq \mc{M}_1(\mc{D}, g)$. This settles local boundedness and Lipschitz regularity of $u$ in $x$ on $\mc{D}$. 
For H\"older regularity, we prepare the estimates required to apply Lemma~\ref{lem:holder_interpolation}.
First, it follows by Banach--Alaoglu in the Bochner space $L^{\infty}(t_0,T;L^2(a,b))$ 
plus the strong convergence of Corollary~\ref{cor:gradients_converge_locally} 
that $\partial_{x}u$ also enjoys the bound of Proposition~\ref{prop:improved_estimates}, that is,
\[
  \operatorname*{ess\,sup}_{t_0<t<T} \normsq{\partial_{x}u(t, \cdot)}_{L^2(a,b)} \leq B_1(t_0,g,u_0).
\]
Next, combining the estimate on $\diff{t}\left(u_n+\frac{1}{n}\right)^{1/2}$ from~\eqref{eq:l2_time_apriori} with the uniform interior
bound on $u_n$ from Proposition~\ref{prop:interior_bounds_n} yields an $L^2$ bound on $\partial_t u_n$,
\[ 
\begin{split} 
  & \normsq{\partial_t u_{n}}_{L^2(\mc{D};dt\times d\nu_{\lambda})} \leq 
 4\norm{u_n+\tfrac1n}_{L^\infty(\mc{D})}\normsq{\partial_t \left(u_n+\frac{1}{n}\right)^{1/2}}_{L^2(\mc{D};dt\times d\nu_{\lambda})} \\
& \leq \frac{\mc{M}_0(\mc{D})C'_T}{t_0^2} + \omega(n).
\end{split}
\]
The last estimate and the strong local convergence of $u_n$ give $\partial_t u_n \wto \partial_t u$ in $L^2(\mathcal{D};dt \times d\nu_{\lambda})$, implying
$\normsq{\partial_t u}_{L^2(\mc{D}; dt \times d\nu_\lambda)} \leq \mc{M}_0(\mc{D})C'_T/t_0^2$.
We transfer this estimate from $L^2(\mc{D};dt \times d\nu_{\lambda})$ to Lebesgue measure at the cost of an additional factor from~\eqref{eq:lambda_to_lebesgue}, giving
\[ \normsq{\partial_t u}_{L^2(\mc{D})} \leq \rho_{\lambda}(b)\mc{M}_0(\mc{D})C'_T/t_0^2. \]
If either $\partial_xu$ or $\partial_tu$ vanishes identically on $\mc{D}$, the claimed H\"older regularity follows directly from one-dimensional Sobolev embedding. Otherwise, the last bounds are exactly those required
to apply Lemma~\ref{lem:holder_interpolation} on $\mc{D}$. This provides a continuous representative 
$\tilde{u} = u\ a.e$ on $\mc{D}$ such that $\tilde{u} \in C^{\frac{1}{4},\frac{1}{2}}(\overline{\mc{D}})$. Explicit 
seminorm bounds in terms of $\partial_{x}u$ and $\partial_t u$ are given in Lemma~\ref{lem:holder_interpolation} for sufficiently small (in time) parabolic cylinders, and the $L^{\infty}(\mc{D})$ 
bound on $u$ has already been established. As $\mc{D}$ was arbitrary, a 
continuous representative on all of $Q_T$ with $\tilde{u} \in C^{\frac{1}{4},\frac{1}{2}}_{loc}(Q_T)$ follows. The bounds are uniform as $t_2\uparrow T$, giving the asserted regularity up to $T$.
\end{proof}

\subsection{Special Solutions}\label{subsection:bc_special_solutions}
The explicit solutions below show that the smoothing theory of this
section is sharp in two respects. First, the B\'enilan--Crandall
inequality~\eqref{eq:bc_for_u} of
Proposition~\ref{prop:benilan_crandall} is saturated, in both branches
of $K_\theta$~\eqref{eq:defn_k_theta}. Second, the spatial regularity
of Theorem~\ref{thm:regularity} cannot be improved within the weak
solution class of Definition~\ref{defn:weak_solution}: we exhibit solutions
that are Lipschitz in $x$ and no better, with jump-discontinuity of the gradient 
at the free boundary $\partial \{u > 0\} \Subset \Omega$. We see also that when $\lambda = 0$, weak solutions need not be globally
bounded.

\medskip\noindent
Fix $\theta_0 \geq 0$ and insert the separable ansatz
$u_\mu(t,x) = G_\mu(t) F_\mu(x)$, with eigenvalue parameter
$\mu > 0$, into~\eqref{eq:classical_form}. We require the reaction to
be constant where the solution is positive (and zero elsewhere), 
so on the positive set we will have $g(x) \defeq \frac{-\theta_0}{\rho_{\lambda}(x)}$. 
This separates the equation, for any $\lambda \geq 0$, into
\begin{align}
  & G_\mu'(t) = -\mu\, G_\mu(t)^2 - \theta_0\, G_\mu(t), \\
  & F_\mu(x)\bigl(\mu + \rho_\lambda(x)\, F_\mu''(x)\bigr) = 0.
\end{align}
The spatial equation does not involve the reaction, and direct
integration gives, using the shifted variable $p = x + \lambda$ when
$\lambda > 0$,
\begin{equation}\label{eq:F_mu_formula}
  F_\mu = \begin{cases}
    \bigl(\mu \log x + c x + \kappa\bigr)_+, & \lambda = 0, \\[1ex]
    \Bigl(\mu\,\arcosh\dfrac{p}{\lambda}
      + p\Bigl(c - \dfrac{\mu}{\lambda}\arcsec\dfrac{p}{\lambda}\Bigr)
      + \kappa\Bigr)_{\!+}, & \lambda > 0
  \end{cases}
\end{equation}
for some constants $c$ and $\kappa$; $u_\mu$ vanishes off the positive
set of $F_\mu$, which is determined by $c$ and $\kappa$. We may
therefore fix the reaction, setting
$g \defeq -\theta_0/\rho_\lambda$ on that positive set and
$g \defeq 0$ elsewhere, so that $\rho_\lambda g \in L^\infty(\Omega)$
and $\theta = \theta_0$ in~\eqref{eq:defn_theta}. The equation for
$G_\mu$ solves as
\begin{equation}\label{eq:G_damped}
  G_\mu(t)
  = \frac{\theta_0}
  {\bigl(\tfrac{\theta_0}{G_\mu(0)} + \mu\bigr) e^{\theta_0 t} - \mu},
\end{equation}
recovering $G_\mu = \bigl(\mu t + 1/G_\mu(0)\bigr)^{-1}$ as
$\theta_0 \to 0$ and decaying at the exponential rate $\theta_0$ when
$\theta_0 > 0$. We impose $G_\mu(0) > 0$ to avoid finite-time blowup,
enforce nonnegativity, and exclude the constant solution. These
solutions therefore saturate~\eqref{eq:bc_for_u} up to a translation
of time: writing
$e^{\theta_0 s} \defeq 1 + \theta_0/\bigl(\mu\, G_\mu(0)\bigr)$ when
$\theta_0 > 0$, and $s \defeq 1/\bigl(\mu\, G_\mu(0)\bigr)$ when
$\theta_0 = 0$, one checks
\begin{equation*}
  \diff{t} u_\mu(t,x)
  = -\bigl(\mu\, G_\mu(t) + \theta_0\bigr)\, u_\mu(t,x)
  = -\frac{u_\mu(t,x)}{K_{\theta_0}(t + s)}.
\end{equation*}
Both branches of~\eqref{eq:defn_k_theta} are therefore attained. It remains to identify the basic
constraints on $c$ and $\kappa$. \textbf{Case $\bm{\lambda = 0}$}.\ For $u_\mu$ to satisfy
Definition~\ref{defn:weak_solution}, we need
$\norm{F_\mu}_{L^p(0, \infty; d\nu_0)} < \infty$ for $p \in \{1, 3\}$,
forcing $c \leq 0$. The parameter $\kappa$ is absorbed by
the change of variables $y := x \exp(\kappa/\mu)$,
$\beta := -c/(\mu \exp(\kappa/\mu)) \geq 0$, which reads
$F_\mu(y) = \mu (\log y - \beta y)_+$. When
$0 < \beta < \frac{1}{e}$, $F_{\mu}$ is supported on the interval
between the two positive solutions of $\log y_{\ast} = \beta
y_{\ast}$, and when $\beta = 0$, it is unbounded and supported on $(1,\infty)$.
\textbf{Case $\bm{\lambda > 0}$}.\ Membership in the weak solution
class requires $c\leq \pi\mu/(2\lambda)$ and $\lambda c+\kappa\leq0$, by the behavior at infinity and the zero trace at $x=0$, respectively.
Nontrivial profiles with free boundary compactly contained in $\Omega$ are obtained by choosing $0<c<\pi\mu/(2\lambda)$ and $\kappa$ so that the positive set is a nonempty compactly contained interval.

\medskip\noindent 
In both cases, we see that $F_\mu$ (and hence $u_{\mu}$) is smooth on its positive set and meets its nontrivial zero
level set in a Lipschitz fashion, with gradient jumping to zero at the
free boundary, matching the regularity of Theorem~\ref{thm:regularity}. Consequently $u_\mu$ is a weak solution
in the sense of Definition~\ref{defn:weak_solution} on all of $Q_T$,
and not only on its positive set: the flux
$\partial_x(u_\mu^2/2) = u_\mu\, \partial_x u_\mu$ extends
continuously by zero across the free boundary, so the weak formulation
follows by integration by parts on the positive set.

\section{Induced Solution of the System}\label{section:system}
This section proves Theorem~\ref{thm:system}. The resonance change of
variables~\eqref{eq:resonance_map} is used to transfer existence and regularity from the scalar theory of
Sections~\ref{section:existence}--\ref{section:regularity} to the
system (Proposition~\ref{prop:regularity_induced_pair},
Corollary~\ref{cor:positive_time_properties}). Under additional 
nonnegativity and moment hypotheses on the initial data, we 
prove that the induced pair is physically admissible, and
Subsection~\ref{subsection:system_examples} exhibits two admissible
solutions. Throughout, $u$ is identified at positive times with its
continuous representative from Theorem~\ref{thm:regularity}, and
where there can be no ambiguity we use the shift of momentum
coordinate $p = x + \lambda$ to write functions of $p$ and $x$
interchangeably.

\subsection{Change of Variables}
Recall the resonance change of variables~\eqref{eq:resonance_map} gives
\begin{equation}\label{eq:inverse_cov}
  x(k) = (k^{-2}+\lambda^{2})^{1/2}-\lambda, \quad
  x'(k) = -\,\frac{k^{-3}}{(k^{-2}+\lambda^{2})^{1/2}},
\end{equation}
a strictly decreasing diffeomorphism of $(0,\infty)$ with inverse
$k(x) = \bigl((x+\lambda)^{2}-\lambda^{2}\bigr)^{-1/2}$, under which
\begin{equation}\label{eq:measure_pullback}
  d\nu_\lambda(x) = \frac{dx}{\rho_\lambda(x)}
  \equiv \frac{dk}{1+\lambda^{2}k^{2}},
  \quad x = x(k).
\end{equation}
We recall the reduction of~\eqref{eq:system} to~\eqref{eq:classical_form} announced in Section~\ref{section:introduction}, 
following~\cite{ivanov1967quasilinear,huang_gamba}; the manipulations below are formal. Beginning
with~\eqref{eq:system}, restrict the momentum coordinate to
$p > \lambda$, where the resonance condition
$(1+\lambda^{2}k^{2})^{1/2} = pk$ has the single root $k = k(p)$,
implying the change of variables for the Dirac delta
in~\eqref{eq:system},
\begin{equation*}
  \delta_0\bigl((1+\lambda^{2}k^{2})^{1/2} - pk\bigr)
 = \frac{p}{p^{2}-\lambda^{2}}\;\delta_0\bigl(k - k(p)\bigr).
\end{equation*}
Writing $u \defeq k^{3}W(k(p),t)$ and
$\gamma(p) \defeq p(p^{2}-\lambda^{2})^{1/2}$, evaluation of the
integrals in~\eqref{eq:system} gives the local pair
\begin{equation}\label{eq:ql_reduced}
  \partial_t f = \partial_p\bigl(\gamma\,u\,\partial_p f\bigr),
  \qquad
  \partial_t u = \gamma\,(\partial_p f)\,u.
\end{equation}
Define $g \defeq \partial_p\bigl(f_0 - \partial_p u_0\bigr)$ and
integrate the above in time to obtain
$\partial_p f = \partial_p^{2} u + g$ and
$f = f_0 + \partial_p(u - u_0)$. Inserting these identities into the
local pair and shifting $x = p-\lambda$
yields~\eqref{eq:classical_form}.

\subsection{The Induced Pair as a Weak Solution}
We now use the induced pair~\eqref{defn:induced_pair} to complete
the converse procedure, wherein we exhibit a weak solution to~\eqref{eq:system} 
starting from the solution $u(t,x)$ of
Theorem~\ref{thm:existence}.

\begin{proposition}[Regularity of the Induced Pair]
\label{prop:regularity_induced_pair}
The induced pair $(f^{\ast},W^{\ast})$ of~\eqref{defn:induced_pair}
satisfies:
\begin{enumerate}
  \item[(i)] for a.e.\ $t \in (0,T)$ the map
  $k \mapsto k^{3}W^{\ast}(t,k)$ lies in
  $L^{q}\bigl((0,\infty); \tfrac{dk}{1+\lambda^{2}k^{2}}\bigr)$ for
  every $q \in [1,3]$, uniformly in $t$, and additionally for
  $q = \delta$ when $\lambda = 0$; moreover $k^{3}W^{\ast}$ admits a
  weak derivative in $k$, square integrable against
  $k^{2}\sqrt{1+\lambda^{2}k^{2}}\,dk\,dt$;
  \item[(ii)] $f^{\ast} - h = \partial_x u\big|_{x = p-\lambda}
  \in L^{2}\bigl((0,T)\times(\lambda,\infty)\bigr)$; in particular
  $f^{\ast}(t,\cdot) - h$ is a function for a.e.\ $t$, even when the
  induced initial distribution $f^{\ast}_0 = h + \partial_p u_0$ is only a
  distribution.
\end{enumerate}
\end{proposition}

\begin{proof}
Inserting~\eqref{eq:measure_pullback} and
$u(t,x(k)) = k^{3}W^{\ast}(t,k)$ from~\eqref{defn:induced_pair}, for
$q \in (0,\infty)$ and a.e.\ $t$,
\begin{equation}\label{eq:spectral_lq}
  \norm{u(t,\cdot)}^{q}_{L^{q}(\Omega;\, d\nu_\lambda)}
  = \int_0^{\infty} k^{3q}\bigl|W^{\ast}(t,k)\bigr|^{q}\,
  \frac{dk}{1+\lambda^{2}k^{2}},
\end{equation}
and, since $x'(k)$ in~\eqref{eq:inverse_cov} carries $dx$ onto
$k^{-2}(1+\lambda^{2}k^{2})^{-1/2}\hspace{0.1em}dk$,
\begin{equation}\label{eq:spectral_grad}
  \norm{\partial_x u}^{2}_{L^{2}(Q_T)}
  = \int_0^T \!\! \int_0^{\infty}
  \bigl|\partial_k\bigl(k^{3}W^{\ast}(t,k)\bigr)\bigr|^{2}\,
  k^{2}\sqrt{1+\lambda^{2}k^{2}}\;dk\,dt.
\end{equation}
The memberships in (i) are then exactly those of
Theorem~\ref{thm:existence}. (ii) is the definition
in~\eqref{defn:induced_pair} under the shift $x = p - \lambda$, with
$\partial_x u \in L^{2}(Q_T)$ again by Theorem~\ref{thm:existence}.
\end{proof}

At positive times the pair inherits more from
Theorem~\ref{thm:regularity} and the B\'enilan--Crandall
inequality~\eqref{eq:bc_for_u}.

\begin{corollary}[Positive-Time Properties]
\label{cor:positive_time_properties}
The induced pair satisfies:
\begin{enumerate}
  \item[(i)] for every $0 < t_0 < T$ and every $K \Subset (0,\infty)$,
  \begin{equation*}
    \norm{W^{\ast}}_{C^{1/4,\,1/2}((t_0,T)\times K)} < \infty
    \quad \text{and} \quad
    \norm{\partial_k W^{\ast}}_{L^{\infty}((t_0,T)\times K)} < \infty,
  \end{equation*}
  and $\norm{f^{\ast} - h}_{L^{\infty}((t_0,T)\times P)} < \infty$ for
  every $P \Subset (\lambda,\infty)$, so $f^{\ast}(t,\cdot)$ is a
  locally bounded function at every positive time;
  \item[(ii)] with $\theta$ and $K_\theta$ as in~\eqref{eq:defn_theta}
  and~\eqref{eq:defn_k_theta},
  $\partial_t W^{\ast} \geq -W^{\ast}/K_\theta(t)$ in the sense of
  distributions on $(0,T)\times(0,\infty)$, and pointwise, for every
  $k \in (0,\infty)$ and $0 < t_1 \leq t_2 \leq T$,
  \begin{equation}\label{eq:spectrum_persistence}
    W^{\ast}(t_2,k) \;\geq\; \sigma_\theta(t_1,t_2)\, W^{\ast}(t_1,k),
    \qquad
    \sigma_\theta(t_1,t_2) \defeq
    \begin{cases}
      t_1/t_2, & \theta = 0,\\[0.5ex]
      \dfrac{e^{\theta t_1}-1}{e^{\theta t_2}-1}, & \theta > 0;
    \end{cases}
  \end{equation}
  in particular the positivity set $\{k : W^{\ast}(t,k) > 0\}$ is
  nondecreasing in $t$ on $(0,T]$, and for $\theta > 0$ the large-time
  damping rate of $W^{\ast}$ is at most $\theta$;
  \item[(iii)] with
  $A(t) \defeq \norm{\partial_x u(t,\cdot)}_{L^{2}(\Omega)}$, which
  satisfies $A\in L^2(0,T)\cap L^\infty_{loc}((0,T])$, for a.e.\
  $t \in (0,T)$, and every $k \in (0,\infty)$,
  \begin{equation*}
    W^{\ast}(t,k) \leq \frac{A(t)}{k^{3}}\,x(k)^{1/2}
    \leq
    \begin{cases}
      (2\lambda)^{-1/2}\,A(t)\,k^{-4}, & \lambda > 0,\\[0.5ex]
      A(t)\,k^{-7/2}, & \lambda = 0.
    \end{cases}
  \end{equation*}
\end{enumerate}
\end{corollary}
\begin{proof}
\emph{(i).} For fixed $K \Subset (0,\infty)$ and $0 < t_0 < T$, the
change of variables $x(k)$ is bi-Lipschitz and carries $K$ onto its
compact image $D = x(K) \Subset \Omega$. On $D$,
Theorem~\ref{thm:regularity} bounds $u$, $\partial_x u$, and the
parabolic H\"older seminorm of $u$, and interior H\"older control of
$W^{\ast}(t,k) = k^{-3}u(t,x(k))$ follows; the Lipschitz bound in $k$
is a composition of Lipschitz maps, with $x'(k)$ as
in~\eqref{eq:inverse_cov}. The same bounds give
$\norm{f^{\ast} - h}_{L^{\infty}((t_0,T)\times P)}
= \norm{\partial_x u}_{L^{\infty}((t_0,T)\times(P-\lambda))} < \infty$.

\emph{(ii).} Testing~\eqref{eq:bc_for_u} with
$\varphi(t,x) = \frac{|k'(x)|}{k(x)^{3}}\,\psi(t,k(x))$, changing
variables $k = k(x)$ and using $u(t,x) = k(x)^{3}W^{\ast}(t,k(x))$,
gives the distributional inequality. For the pointwise bound, fix
$t_0 \in (0,T)$ and set
$E(t) \defeq \exp\int_{t_0}^{t} K_\theta(s)^{-1}ds$, so
$E' = E/K_\theta$ and $\langle \partial_t(uE), \psi\rangle
= \langle \partial_t u + u/K_\theta,\, E\psi\rangle \geq 0$
by~\eqref{eq:bc_for_u}. So $uE$ is nondecreasing in $t$ (by
mollification in $x$ and the continuity of $u$ from
Theorem~\ref{thm:regularity}), whence
$u(t_2,x) \geq u(t_1,x)E(t_1)/E(t_2)$ with
$E(t_1)/E(t_2) = \sigma_\theta(t_1,t_2)$, since
$K_0(s)^{-1} = 1/s$ and
$K_\theta(s)^{-1} = \tfrac{d}{ds}\log(e^{\theta s}-1)$ for
$\theta > 0$. Multiplying by $k^{-3}$ at $x = x(k)$
gives~\eqref{eq:spectrum_persistence}; the positivity set is then
nondecreasing, and
$\sigma_\theta \geq (1-e^{-\theta t_1})e^{-\theta(t_2-t_1)}$ gives the
exponential rate.

\emph{(iii).} As in the proof of Theorem~\ref{thm:regularity}, the
estimates of Proposition~\ref{prop:improved_estimates} transfer to $u$
in the limit, giving $\operatorname*{ess\,sup}_{t\in(t_0,T)}A(t)^2\leq B_1(t_0,g,u_0)$ for every $t_0>0$, while $A\in L^2(0,T)$ follows from~\eqref{eq:lp_apriori}. Since
$u(t,\cdot) \in V_0(\Omega)$ has vanishing trace at $x = 0$, the
fundamental theorem of calculus and Cauchy--Schwarz give
$u(t,x) \leq A(t)\,x^{1/2}$ for a.e.\ $t$, and applying the resonance
change of variables,
\[
  W^{\ast}(t,k) \;\leq\; \frac{A(t)}{k^{3}}\,
  \Bigl( (\lambda^{2}+k^{-2})^{1/2}-\lambda \Bigr)^{1/2}.
\]
For $\lambda = 0$ the inner factor is $k^{-1}$ exactly; for
$\lambda > 0$, concavity of the square root gives
$(\lambda^{2}+k^{-2})^{1/2}-\lambda \leq k^{-2}/(2\lambda)$.
\end{proof}

Note that by~\eqref{eq:spectral_lq}, the data class of
Theorem~\ref{thm:existence} reads
$k^{3}W^{\ast}_0 \in L^{1} \cap
L^{3}\bigl(\tfrac{dk}{1+\lambda^{2}k^{2}}\bigr)$ for $\lambda > 0$,
and $k^{3}W^{\ast}_0 \in L^{\delta}(dk) \cap L^{3}(dk)$ for $\lambda = 0$.
With the change of variables absorbing the measure $d\nu_\lambda$, the
notion of weak solution from Definition~\ref{defn:weak_solution}
transfers to the system. First, we define
\begin{equation}\label{defn:gamma}
  \gamma(p) \defeq \rho_\lambda(p-\lambda) = p(p^{2}-\lambda^{2})^{1/2}
\end{equation}
and the interior test classes
\begin{equation*}
  \mathcal{Z} \defeq C^{\infty}_c\bigl((\lambda,\infty) \times [0,T)\bigr),
  \qquad
  \Xi \defeq C^{\infty}_c\bigl((0,\infty) \times [0,T)\bigr).
\end{equation*}

\begin{proposition}[Weak Solution of the Quasilinear System]
\label{prop:system_weak}
With $u_0$, $g$, and $u$ as in Theorem~\ref{thm:existence}, and
$(f^{\ast},W^{\ast})$ given by~\eqref{defn:induced_pair}, the
identities~\eqref{eq:system_weak_intro} hold for every
$\zeta \in \mathcal{Z}$ and $\xi \in \Xi$, with the $f_0^{\ast}$ term interpreted by distributional duality and all other integrals absolutely convergent.
\end{proposition}

\begin{proof}
By~\eqref{defn:induced_pair}, $u(t,p-\lambda) =
(p^{2}-\lambda^{2})^{-3/2}\,W^{\ast}(t,k(p))$. Hence
\begin{equation}\label{eq:kernel_weight}
  \gamma(p)\,u(t,p-\lambda)
  = \frac{p}{p^{2}-\lambda^{2}}\;W^{\ast}(t,k(p)),
\end{equation}
so that the right-hand sides of~\eqref{eq:system_weak_intro} may
equivalently be written with the kernel
$\tfrac{p}{p^{2}-\lambda^{2}}\,W^{\ast}(t,k(p))$ in place of
$\gamma\,u$, the form of~\cite[Eqs.~(3.1), (3.3)]{huang_gamba}; we
work with $u$ throughout. Let
$\eta \in C^{\infty}_c([0,T)\times\Omega)$. Since $\rho_\lambda$ is
smooth and positive on the open half-line, the product
$\rho_\lambda\eta$ is admissible in
Definition~\ref{defn:weak_solution}. Inserting it
into~\eqref{eq:weak_form} and expanding
$\partial_x(\rho_\lambda\eta)$ gives
\begin{equation}\label{eq:master_identity}
  \begin{split}
  &-\int_{Q_T} u\,\partial_t\eta
  + \int_{Q_T} \frac{\rho_\lambda}{2}\,\partial_x u^{2}\,\partial_x\eta
  + \int_{Q_T} \frac{\partial_x\rho_\lambda}{2}\,\partial_x u^{2}\,\eta
  + \int_{Q_T} \rho_\lambda\,|\partial_x u|^{2}\,\eta \\
  &\qquad = \int_{\Omega} u_0\,\eta(0,\cdot)
  + \int_{Q_T} \rho_\lambda g\, u\,\eta,
  \end{split}
\end{equation}
with all integrals finite by Theorem~\ref{thm:existence}. To
prove the first identity of~\eqref{eq:system_weak_intro}, take
$\eta \defeq \partial_p\zeta$, which
lies in $C^{\infty}_c([0,T)\times\Omega)$ because
$\support \zeta$ is compact in $(\lambda,\infty)$.
Substitute $f^{\ast} = h + \partial_p u$ and
$f^{\ast}_0 = h + \partial_p u_0$ on its left-hand side. The
$h$-contributions cancel, because $h$ is time-independent and
$\zeta(T,\cdot) = 0$. Integrating by parts in $p$, that left-hand side
becomes
$\int_{Q_T} u\,\partial_t\eta + \int_{\Omega} u_0\,\eta(0,\cdot)$. On
the right, the same substitution, the identity
$u\,\partial_p u = \tfrac12 \partial_p u^{2}$, and one integration by
parts in the $h$ term give
\begin{equation*}
  \int_{Q_T} \Bigl(
  \frac{\rho_\lambda}{2}\,\partial_x u^{2}\,\partial_x\eta
  + \frac{\partial_x\rho_\lambda}{2}\,\partial_x u^{2}\,\eta
  + \rho_\lambda\,|\partial_x u|^{2}\,\eta
  - \rho_\lambda g\, u\,\eta \Bigr).
\end{equation*}
This integration by parts is valid because
$\partial_p h = g(\cdot-\lambda)$ and $\gamma u\eta$ belongs to $H^1$ in space for a.e.\ $t$, with compact spatial support in $\Omega$ and the required space--time integrability. Thus the first identity
of~\eqref{eq:system_weak_intro} is precisely~\eqref{eq:master_identity}
rearranged.

For the second identity, take $\eta_\xi(t,p) \defeq p\,\xi(k(p),t)$.
This again lies in $C^{\infty}_c([0,T)\times\Omega)$, because the
resonance map carries $\support \xi \Subset (0,\infty)$ onto
a compact subset of $(\lambda,\infty)$. Under the change of variables
$k = k(p)$, the Jacobian $|dk/dp| = p\,(p^{2}-\lambda^{2})^{-3/2}$
cancels the factor
$W^{\ast}(t,k(p)) = (p^{2}-\lambda^{2})^{3/2}u(t,p-\lambda)$ and
leaves the weight $p$. Hence
\begin{equation*}
  \int_0^{\infty} W^{\ast}\,\partial_t\xi\;dk
  = \int_\lambda^{\infty} u\,\partial_t\eta_\xi\;dp,
  \qquad
  \int_0^{\infty} W^{\ast}_0\,\xi(\cdot,0)\;dk
  = \int_\lambda^{\infty} u_0\,\eta_\xi(0,\cdot)\;dp.
\end{equation*}
Applying the argument for the first identity with $\eta$ replaced by
$\eta_\xi$ completes the proof.
\end{proof}

\subsection{Physical Admissibility}
Physical admissibility of the induced pair demands $W^{\ast} \geq 0$,
$f^{\ast}(t,\cdot) \geq 0$, and unit mass
$\int_\lambda^\infty f^{\ast}(t,p)\,dp = 1$. While nonnegativity of $W^{\ast}$ is immediate from $u \geq 0$
(Theorem~\ref{thm:existence}), proving that admissibility of
$f^{\ast}$ propagates from the initial data to positive times is the
content of the next two propositions.
\begin{proposition}[Nonnegativity of Induced $f^{\ast}$]
\label{prop:f_nonneg}
Let $u_0$ and $g$ satisfy the hypotheses of Theorem~\ref{thm:existence},
with $h$ as in~\eqref{eq:coupling_profile}. Writing
$f^{\ast}_0 \defeq h + \partial_p u_0$ as consistent
with~\eqref{defn:induced_pair}, suppose also that $f^{\ast}_0 \geq 0$ in
$\mathcal{D}'\bigl((\lambda,\infty)\bigr)$, as
in~\eqref{eq:f0_positive_hypothesis}. Then one may choose particular regularized data $u_{0,n}$, $g_n$
compatible with Definition~\ref{defn:n_problem} such that the
constructed solution $u$ of Theorem~\ref{thm:existence} gives
\begin{equation}\label{eq:f_positive_conclusion}
  f^{\ast}(t,x) = h(x) + \diff{x} u(t,x) \;\geq\; 0
  \quad \text{a.e.\ on } Q_T.
\end{equation}
\end{proposition}

\begin{proof}
Denote the negative part $(s)_- \defeq \max(-s, 0)$, and for data
$(u_{0,n}, g_n)$ admissible in Definition~\ref{defn:n_problem}, with
$u_n$ as in Proposition~\ref{prop:n_are_classical}, write
$h_n \defeq -\int_x^\infty g_n$, $f_{0,n} \defeq h_n + \diff{x}u_{0,n}$,
and $f_n \defeq h_n + \diff{x} u_n$. It suffices to prove, for every
fixed $K \Subset Q_T$, that $\norm{(f_n)_-}_{L^2(K)} \to 0$: by
Corollary~\ref{cor:gradients_converge_locally} and the locally
uniform convergence $h_n \to h$ established below,
$(f_n)_- \to (f^\ast)_- = (h + \diff{x} u)_-$ in $L^2(K)$ along the
(not relabeled) subsequence of the construction.

To this end we first show one may choose $f_{0,n}$ asymptotically
nonnegative, i.e.\ $\norm{(f_{0,n})_-}_{L^2(\Omega_n)} = \omega(n)$.
To see this, take $r_n=n^{-\delta/2}$ as allowed in Definition~\ref{defn:n_problem} and
$(a_n, b_n) \defeq (r_n^{1/4},\, 3n)$, admissible
in~\eqref{eq:omega_n_requirement}. Define also
$\chi_n \in C^\infty(\Omega; [0,1])$ with $\chi_n \equiv 1$ on
$[a_n + 3r_n,\, n]$, $\chi_n \equiv 0$ on
$(0,\infty) \setminus [a_n + 2r_n,\, 2n]$, with
$\abs{\chi_n'} \leq 2/r_n$ on the ascent and
$\abs{\chi_n'} \leq 2/n$ on the descent $[n, 2n]$. Take a sequence of
mollification scales $\varepsilon_n \in (0, r_n]$ sufficiently small, using a nonnegative mollifier and extension by zero, and now define
\begin{equation}\label{eq:data_defn}
  u_{0,n} \defeq (\chi_n u_0)_{\varepsilon_n},
  \quad
  g_n \defeq \rho_\lambda^{-1}(\chi_n\, \rho_\lambda g)_{\varepsilon_n},
\end{equation}
so that $\support u_{0,n} \cup \support g_n
\subset [a_n + r_n,\, b_n - r_n]$. These data satisfy~\eqref{eq:n_problem_data_requirement}; in particular convolution preserves $\rho_\lambda g_n\geq-\theta$ and $\linfty{\rho_\lambda g_n}\leq\linfty{\rho_\lambda g}$. Note
\begin{equation}\label{eq:descent_L2}
  \norm{u_0 \chi_n'}_{L^2(n, 2n)}
  \leq \frac{2}{n}\, \norm{u_0}_{L^2(n, 2n)}
  \leq C_\lambda
  \Bigl(\int_n^\infty u_0^3\, d\nu_\lambda\Bigr)^{\!1/3} n^{-1/6}
  = \omega(n),
\end{equation}
and that
\begin{equation}\label{eq:h_tail}
  \abs{h(x)}
  \leq \linfty{\rho_\lambda g}\, \nu_\lambda\bigl((x, \infty)\bigr)
  \leq C_\lambda \linfty{\rho_\lambda g}\, x^{-1},
  \qquad x > 0,
\end{equation}
valid also for $h_n$, since
$\linfty{\rho_\lambda g_n} \leq \linfty{\rho_\lambda g}$. By the $BV$
product rule ($Du_0 = f^{\ast}_0 - h\,dx$ as measures, $u_0 \in BV_{loc}$
by~\eqref{eq:f0_positive_hypothesis} and~\eqref{eq:h_tail}) and
mollification,
\begin{equation*}
  f_{0,n}
  = (\chi_n f^{\ast}_0)_{\varepsilon_n}
  + (u_0 \chi_n')_{\varepsilon_n}
  + \bigl[h_n - (\chi_n h)_{\varepsilon_n}\bigr].
\end{equation*}
The first term is nonnegative due to the hypothesis on $f^{\ast}_0$; the
second is nonnegative outside the $\varepsilon_n$-neighborhood
$[n - \varepsilon_n,\, 2n + \varepsilon_n]$ of the descent, on which
it is $\omega(n)$ in $L^2$ by~\eqref{eq:descent_L2}; and
by~\eqref{eq:h_tail} the third is $\omega(n)$ in $L^2(\Omega_n)$ for
$\varepsilon_n$ small. Indeed, before mollification the ascent contributes at most $C\sqrt{r_n}/a_n=C r_n^{1/4}$ and the upper cutoff at most $C n^{-1/2}$; the mollification errors can be made $\omega(n)$. The same bounds give $h_n \to h$ locally uniformly on $\Omega$.

To conclude, we now use the regularity of $u_n$ from
Proposition~\ref{prop:n_are_classical}. Since
$\diff{x} f_n = \diff{x}^2 u_n + g_n$, we find
\begin{equation}\label{eq:flux_form}
  \diff{t} u_n
  = \rho_\lambda \bigl(u_n + \tfrac{1}{n}\bigr)\, \diff{x} f_n
  \qquad \text{on } \overline{Q_{T,n}}.
\end{equation}
As $\diff{t} f_n$ may not exist pointwise, test~\eqref{eq:flux_form} with
$\varphi \in H^1(\Omega_n)$,
$\eta \in C^\infty_c\bigl((0,T)\bigr)$ to find
\begin{equation}\label{eq:fn_weak_form}
  \begin{split}
  -\int_0^T \!\eta' \int_{\Omega_n} f_n \varphi\,dx\,dt
  &= \int_0^T \!\eta' \int_{\Omega_n} u_n\, \diff{x}\varphi\,dx\,dt \\
  &= -\int_0^T \!\eta \int_{\Omega_n}
    \rho_\lambda \bigl(u_n + \tfrac{1}{n}\bigr)\, \diff{x} f_n\,
    \diff{x}\varphi\,dx\,dt.
  \end{split}
\end{equation}
We conclude
$\diff{t} f_n \in L^\infty\bigl(0,T; H^1(\Omega_n)^*\bigr)$, i.e.\
against all of $H^1(\Omega_n)$, and since
$f_n \in C([0,T]; L^2(\Omega_n)) \cap L^2(0,T; H^1(\Omega_n))$, we
find $\varphi = -(f_n)_-(t)$ is a valid test function. The weak chain
rule (with $s \mapsto \tfrac{1}{2}s_-^2$, see, e.g.,
\cite[Lemma~1.5]{AltLuckhaus1983}) then gives
$\tfrac{d}{dt}\tfrac{1}{2}\norm{(f_n)_-(t,\cdot)}^2_{L^2(\Omega_n)}
= -\int_{\Omega_n} \rho_\lambda \bigl(u_n + \tfrac{1}{n}\bigr)
\abssq{\diff{x}(f_n)_-} \leq 0$, using $u_n+1/n>0$ from~\eqref{eq:shifted_positivity}, and after integration,
\begin{equation}\label{eq:fn_contraction}
  \norm{(f_n)_-(t,\cdot)}_{L^2(\Omega_n)}
  \leq \norm{(f_{0,n})_-}_{L^2(\Omega_n)},
  \qquad t \in [0,T].
\end{equation}
As $\norm{(f_{0,n})_-}_{L^2(\Omega_n)} = \omega(n)$ and
$\norm{(f_n)_-}^2_{L^2(K)}
\leq T \sup_t \norm{(f_n)_-(t,\cdot)}^2_{L^2(\Omega_n)}$, the proof
is complete.
\end{proof}

\begin{proposition}[Conservation of Mass of $f^{\ast}$]\label{prop:f_mass}
In addition to the hypotheses of Proposition~\ref{prop:f_nonneg},
suppose the moment
conditions~\eqref{eq:moment_hypotheses} hold, so that
$h \in L^1\bigl((\lambda,\infty)\bigr)$, and suppose further the unit
normalization~\eqref{eq:unit_mass}. Then the data of Definition~\ref{defn:n_problem} may be chosen so
that the constructed solution $u$ of Theorem~\ref{thm:existence}
satisfies, in addition to~\eqref{eq:f_positive_conclusion},
\begin{equation}\label{eq:f_mass}
  \int_\lambda^\infty f^{\ast}(t, p)\,dp = 1
  \qquad \text{for a.e. } t \in (0, T).
\end{equation}
\end{proposition}
\begin{proof}
Take $h_n \defeq \kappa_n - \int_x^\infty g_n$, with $\kappa_n$ enforcing
$\int_{\Omega_n} h_n\,dx = 1$ (valid as $\Omega_n$ has finite
measure). This constant shift leaves~\eqref{eq:fn_weak_form} unchanged; its contribution to the initial negative part is estimated below.
Since $u_n$ vanishes on $\partial\Omega_n$ and
$f_n = h_n + \partial_x u_n$,
\begin{equation*}
  \int_{\Omega_n} f_n(t,\cdot)\,dx
  = \int_{\Omega_n} h_n\,dx
  + \bigl[u_n(t,\cdot)\bigr]_{a_n}^{b_n}
  = 1
  \qquad \text{for every } n \text{ and every } t,
\end{equation*}
so mass is conserved exactly for the regularized problems. It remains
to pass to the limit without losing mass at either end of $\Omega$.

\smallskip \noindent 
To do this, we first show the negative part of $f_n$ vanishes as $n \to \infty$. Let $\Psi_\varepsilon \in C^{1,1}(\bbR)$ be convex with
$\Psi_\varepsilon'$ Lipschitz, $\Psi_\varepsilon(s) \to s_-$ locally
uniformly as $\varepsilon \to 0$, and $\Psi_\varepsilon'' \geq 0$
a.e. Testing~\eqref{eq:fn_weak_form} with $\Psi_\varepsilon'(f_n)$,
the weak chain rule~\cite[Lemma~1.5]{AltLuckhaus1983} shows
$t \mapsto \int_{\Omega_n}\Psi_\varepsilon(f_n(t,\cdot))$ is
absolutely continuous with
\begin{equation*}
  \frac{d}{dt}\int_{\Omega_n}\Psi_\varepsilon(f_n)
  = -\int_{\Omega_n} \rho_\lambda\bigl(u_n+\tfrac1n\bigr)\,
    \Psi_\varepsilon''(f_n)\,\abssq{\partial_x f_n} \leq 0
  \qquad \text{for a.e. } t.
\end{equation*}
The sign again follows from~\eqref{eq:shifted_positivity}. Integrating in time and then taking $\varepsilon \to 0$ (the latter by
dominated convergence) gives the $L^1$ analogue
of~\eqref{eq:fn_contraction},
\begin{equation}\label{eq:fn_contraction_L1}
  \norm{(f_n)_-(t,\cdot)}_{L^1(\Omega_n)}
  \leq \norm{(f_{0,n})_-}_{L^1(\Omega_n)}
  \qquad \text{for every } t \in [0,T].
\end{equation}
Now we show $\norm{(f_{0,n})_-}_{L^1(\Omega_n)} = \omega(n)$ for an
appropriate choice of the regularized data.  As in the
proof of Proposition~\ref{prop:f_nonneg}, with $\chi_n$,
$\varepsilon_n$, and the data~\eqref{eq:data_defn} as chosen there, we
decompose
\begin{equation*}
  f_{0,n}
  = \bigl(\chi_n f^{\ast}_0\bigr)_{\varepsilon_n}
  + \bigl(u_0 \chi_n'\bigr)_{\varepsilon_n}
  + \bigl[h_n - (\chi_n h)_{\varepsilon_n}\bigr],
\end{equation*}
where the first term is nonnegative by hypothesis~\eqref{eq:f0_positive_hypothesis}. Choose $\varepsilon_n$ still smaller if necessary so that $\int_0^\infty y\,|g_n-\chi_n g|\,dy=\omega(n)$; dominated convergence and~\eqref{eq:moment_hypotheses} then give $\int_0^\infty y\,|g-g_n|\,dy=\omega(n)$. The third term is
$\omega(n)$ in $L^1(\Omega_n)$: by Fubini
and~\eqref{eq:moment_hypotheses},
$\abs{\kappa_n}\abs{\Omega_n}
\leq \int_0^\infty y\,\abs{g - g_n}\,dy
+ \int_{\Omega\setminus\Omega_n}\abs{h}\,dx = \omega(n)$, and consequently
$\norm{h_n-h}_{L^1(\Omega_n)}\leq |\kappa_n||\Omega_n|+\int_0^\infty y\,|g-g_n|\,dy=\omega(n)$, while
$\norm{h - (\chi_n h)_{\varepsilon_n}}_{L^1(\Omega_n)} = \omega(n)$
by $h \in L^1$, dominated convergence, and $\varepsilon_n$ small.

\smallskip\noindent
It remains to estimate $u_0(\chi_n')_-$, supported on
$[n, 2n]$, where $\abs{\chi_n'} \leq 2/n$. By H\"older's inequality
with exponents $(3, 3/2)$,
\begin{equation*}
  \norm{u_0(\chi_n')_-}_{L^1}
  \leq \frac{2}{n} \int_n^{2n} u_0\,dy
  \leq \frac{2}{n}\, n^{2/3}
    \Bigl(\int_n^{2n} u_0^3\,dy\Bigr)^{1/3}
  = 2\Bigl(\frac{1}{n}\int_n^{2n} u_0^3\,dy\Bigr)^{1/3}.
\end{equation*}
On $[n,2n]$ one has $\rho_\lambda(y) \leq C_\lambda\, n\, y$, so that 
\begin{equation*}
  \frac{1}{n}\int_n^{2n} u_0^3\,dy
  \leq C_\lambda \int_n^{2n} y\, u_0^3\, d\nu_\lambda(y)
  \leq C_\lambda \int_n^\infty y\, u_0^3\, d\nu_\lambda(y)
  = \omega(n),
\end{equation*}
by the second moment condition of~\eqref{eq:moment_hypotheses}.
Nonnegative mollification gives
\[
 \norm{\bigl((u_0\chi_n')_{\varepsilon_n}\bigr)_-}_{L^1}
 \leq\norm{u_0(\chi_n')_-}_{L^1}=\omega(n).
\]
Altogether $\norm{(f_{0,n})_-}_{L^1(\Omega_n)} = \omega(n)$. Combining
with~\eqref{eq:fn_contraction_L1}, we find
\[
 \norm{(f_n)_-(t,\cdot)}_{L^1(\Omega_n)} = \omega(n)
 \qquad\text{uniformly in }t.
\]

\smallskip \noindent
Now we can pass to the limit $n \to \infty$. Fix $R > 1$ and set $K_R \defeq [1/R,\, R]$. For the upper bound,
\begin{equation*}
  \int_{K_R} f_n(t,\cdot)\,dx
  = 1 - \int_{\Omega_n \setminus K_R} f_n(t,\cdot)\,dx
  \leq 1 + \norm{(f_n)_-(t,\cdot)}_{L^1(\Omega_n)}
  = 1 + \omega(n).
\end{equation*}
For the lower bound, the derivative part of the tail telescopes,
$\int_{\Omega_n\setminus K_R}\partial_x u_n
= u_n(t,1/R) - u_n(t,R)$, and using $u_n(t,R)\geq-1/n$ from~\eqref{eq:n_maximum_principle},
\begin{equation*}
  \int_{K_R} f_n(t,\cdot)\,dx
  \geq 1 - u_n(t,1/R) - \omega(n)
  - \int_{\Omega\setminus K_R}\abs{h}\,dx,
\end{equation*}
where the $\omega(n)$ collects the $1/n$ error and the $h_n$-tail via
$\norm{h_n - h}_{L^1(\Omega_n)} = \omega(n)$. Since
$u_n(t,\cdot)$ vanishes at $a_n < 1/R$,
Cauchy--Schwarz gives
\[
 u_n(t,1/R) \leq R^{-1/2}\,\norm{\partial_xu_n(t,\cdot)}_{L^2(\Omega_n)},
\]
whose integral in time is bounded
uniformly in $n$ by~\eqref{eq:lp_apriori}. Integrating the two-sided bound over $(s,s') \subset (0,T)$ and
using Corollary~\ref{cor:gradients_converge_locally} (whose local
$L^2$ convergence implies local $L^1$ convergence, as
$(s,s')\times K_R$ has finite measure), together with
$\norm{h_n - h}_{L^1(\Omega_n)} = \omega(n)$, we find
$f_n \to f^{\ast}$ in $L^1\bigl((s,s')\times K_R\bigr)$, implying
\begin{equation*}
  (s'-s)\Bigl(1 - \int_{\Omega\setminus K_R}\abs{h}\,dx\Bigr)
  - C\,R^{-1/2}
  \;\leq\; \int_s^{s'}\!\!\int_{K_R} f^{\ast}\,dx\,dt
  \;\leq\; s'-s.
\end{equation*}
Letting $R \to \infty$, the $h$-tail vanishes since $h \in L^1$, and
$\int_{K_R} f^{\ast} \to \int_\Omega f^{\ast}$ monotonically, by monotone
convergence, using $f^{\ast} \geq 0$ from~\eqref{eq:f_positive_conclusion}. Hence
\begin{equation*}
  \int_s^{s'}\!\!\int_\Omega f^{\ast}\,dx\,dt = s' - s
  \qquad \text{for all } 0 < s < s' < T,
\end{equation*}
and~\eqref{eq:f_mass} follows for a.e.\ $t$ by Lebesgue
differentiation.
\end{proof}

\subsection{Proof of Theorem~\ref{thm:system}}
We are now in position to prove Theorem~\ref{thm:system}, which we
restate here.

\newtheorem*{thmsystemrestated}{Theorem~\ref{thm:system} (Induced
Solution of the System, restated)}
\begin{thmsystemrestated}
Let $u_0$ and $g$ satisfy the hypotheses of
Theorem~\ref{thm:existence}. Then:
\begin{enumerate}
  \item[(i)] the induced pair $(f^{\ast}, W^{\ast})$
  of~\eqref{defn:induced_pair} is a weak solution of the
  system~\eqref{eq:system} in the sense
  of~\eqref{eq:system_weak_intro}, with $W^{\ast} \geq 0$;
  \item[(ii)] if moreover~\eqref{eq:f0_positive_hypothesis} holds,
  then $u$ may be constructed so that additionally
  $f^{\ast}(t,\cdot) \geq 0$ a.e.\ on $(\lambda,\infty)$ for a.e.\
  $t \in (0,T)$;
  \item[(iii)] if additionally~\eqref{eq:moment_hypotheses}
  and~\eqref{eq:unit_mass} hold, the same construction conserves
  mass, $\int_\lambda^\infty f^{\ast}(t,p)\,dp = 1$ for a.e.\
  $t \in (0,T)$.
\end{enumerate}
\end{thmsystemrestated}

\begin{proof}
\emph{(i)} is Propositions~\ref{prop:regularity_induced_pair}
and~\ref{prop:system_weak}, with $W^{\ast} \geq 0$ following from
$u \geq 0$ (Theorem~\ref{thm:existence})
and~\eqref{defn:induced_pair}. \emph{(ii)} is
Proposition~\ref{prop:f_nonneg}, and \emph{(iii)} is
Proposition~\ref{prop:f_mass}, whose data choice refines that of
Proposition~\ref{prop:f_nonneg}.
\end{proof}

\subsection{Admissible Solutions}\label{subsection:system_examples}

Two admissible solutions exhibit the scope of Theorem~\ref{thm:system} in its hypotheses and conclusions. The first
has cold-beam data, an initial distribution $f_0^{\ast}$ which is a measure and
for which no explicit solution is available. The second is explicit
in both variables, conserving mass exactly at every time and
attaining the exponential rate in the persistence bound~\eqref{eq:spectrum_persistence}.

\medskip\noindent
\textbf{Two cold beams.}\ Fix $0 < a < b < c$, let $h \geq 0$ be
smooth, supported in $[c, c+1]$ with $\int_{0}^{\infty} h = 1$, and set
$g \defeq h'$, so that $\rho_\lambda g \in L^\infty$
and~\eqref{eq:moment_hypotheses}--\eqref{eq:unit_mass} hold. Take
initial data
\begin{equation}\label{eq:example_beam_data}
  u_0(x) \defeq \tfrac12\,\mathbf{1}_{[a,\infty)}(x)
  + \tfrac12\,\mathbf{1}_{[b,\infty)}(x)
  - \int_0^x h(y)\,dy.
\end{equation}
Since $\int_{0}^{\infty} h = 1$ is the combined weight of the two jumps, $u_0$ is
the bounded staircase taking the value $\tfrac12$ on $[a,b)$ and $1$
on $[b,c)$, decreasing to $0$ at $x = c+1$. In particular $u_0 \geq 0$
and is supported in $[a,\, c+1]$. The hypotheses of
Theorem~\ref{thm:existence} hold for all $\lambda \geq 0$.
We find induced initial distribution $f^{\ast}_0$, on momentum
$p = x + \lambda$, given by the pair of cold beams
\begin{equation}\label{eq:example_beam_f0}
  f^{\ast}_0
  = \tfrac12\bigl(\delta_{\lambda+a} + \delta_{\lambda+b}\bigr),
\end{equation}
an elementary bump-on-tail configuration of quasilinear
theory~\cite{drummond,vedenov}.
No explicit solution is available here, separability fixing a
continuous profile. The construction nonetheless describes the
evolution completely. By Theorem~\ref{thm:system}, $f^{\ast}(t,\cdot) \geq 0$ has unit mass for a.e.\ $t \in (0,T)$, and by
Corollary~\ref{cor:positive_time_properties}(i) it is a locally
bounded function at every positive time: the two Dirac masses are
replaced instantaneously by a density of the same total mass. Parts~(ii) and~(iii) of Corollary~\ref{cor:positive_time_properties} then
apply, giving persistence of the excited spectrum and its decay at
large wavenumber, the damping rate being at most
$\theta = \esssup(-\rho_\lambda h') > 0$ from~\eqref{eq:defn_theta}, positive since $h$ is nonconstant.

\medskip\noindent
\textbf{An explicit solution.}\ Take $F$ as
in~\eqref{eq:F_mu_formula} with $\mu = 1$: for $\lambda=0$ choose $0<\beta<1/e$, and for $\lambda>0$ choose $0<c<\pi/(2\lambda)$ and
$-\arcosh(\sec(\lambda c))<\kappa<-\lambda c$,
\begin{equation}\label{eq:example_profile}
  F = \begin{cases}
    (\log x - \beta x)_+, & \lambda = 0, \\[1ex]
    \Bigl(\arcosh\dfrac{p}{\lambda}
      + p\Bigl(c - \dfrac{1}{\lambda}\arcsec\dfrac{p}{\lambda}\Bigr)
      + \kappa\Bigr)_{\!+}, & \lambda > 0,
  \end{cases}
\end{equation}
so that $\support F = [\sigma_-, \sigma_+] \Subset \Omega$ with
$F(\sigma_\pm) = 0$. Choose $\theta_0>0$ so that $\theta_0\int_{\sigma_-}^{\sigma_+}y\,\rho_\lambda(y)^{-1}\,dy<1$. With $u(t,x) \defeq G(t) F(x)$ and $G$ as
in~\eqref{eq:G_damped}, set
\begin{equation}\label{eq:example_gh}
  g \defeq -\theta_0\, \rho_\lambda^{-1}\,
  \mathbf{1}_{[\sigma_-,\, \sigma^{\ast}]},
  \qquad
  h(p) = \theta_0 \left(\int_{\max(p - \lambda,\, \sigma_-)}^{\sigma^{\ast}}
  \frac{dy}{\rho_\lambda(y)}\right)_+,
\end{equation}
where $\sigma^{\ast} > \sigma_+$ is the unique value for which
$\theta_0 \int_{\sigma_-}^{\sigma^{\ast}} y\,\rho_\lambda(y)^{-1}\,dy = 1$. Then
$\rho_\lambda g \equiv -\theta_0$ on $\support F$, so $u$ is a damped
separable solution of
Subsection~\ref{subsection:bc_special_solutions}. By Fubini
we see that $\int_\lambda^\infty h\,dp = 1$, the
moments~\eqref{eq:moment_hypotheses} are finite, and
$f^{\ast}_0 = h + G(0) F' \geq 0$ for $G(0)$ small, $h$ being bounded
below and $F'$ bounded on $\support F$. The induced distribution
$f^{\ast}(t,\cdot) = h + G(t) F'$ therefore carries jump
discontinuities at $p = \lambda+\sigma_\pm$ and remains nonnegative for every
$t$, since $G$ decreases. Mass is conserved exactly:
$F(\sigma_\pm) = 0$ gives
$\int_\lambda^\infty f^{\ast}(t,\cdot)\,dp = \int_\lambda^\infty
h\,dp = 1$ for every $t$, the boundary mechanism of
Proposition~\ref{prop:f_mass}. Finally
$\esssup(-\rho_\lambda g) = \theta_0$, so
$\norm{f^{\ast}(t,\cdot) - h}_{L^\infty} \lesssim e^{-\theta_0 t}$
and $W^{\ast}(t,\cdot)$ damps at exactly the rate $\theta_0$, attaining
the exponential rate in Corollary~\ref{cor:positive_time_properties}(ii).

\section*{Appendix} 
\appendix
\renewcommand{\thesection}{\Alph{section}}
\section{Hardy Inequalities}\label{appendix:hardy_sobolev}
We recall a sharp characterization of one-dimensional Hardy inequalities from 
\cite{opichardy}, verifying their particular forms for the weighted Sobolev spaces $V^p((0,b);d\sigma)$, $V^p_L((0, b); d\sigma)$,
and $V^p_0((0, b); d\sigma)$, defined in Subsection~\ref{subsection:notation}. 

\begin{theorem}[1-Dimensional Hardy Inequalities, {\cite[Theorems 5.9, 7.3]{opichardy}}]\label{thm:hardy_opic}
Let $b \leq \infty$ and let $d\sigma = s(x) \, dx$, $d\mu = m(x) \, dx$ be Borel 
measures on $(0, b)$ with $s, m$ finite Lebesgue-a.e. For $1 \leq p \leq q < \infty$ 
and $\frac{1}{p} + \frac{1}{p'} = 1$, the inequality
\begin{equation}\label{eq:hardy_inequality}
  \norm{f}_{L^q(0, b; d\sigma)} \leq C_L(p, q, b) \norm{\partial_{x}f}_{L^p(0, b; d\mu)}
\end{equation}
holds for all
\[
  f \in \mc{D}_L(0,b) = \{f \in C^{\infty}(0,b) \colon \lim_{x\to 0^+} f(x) =0\}
\] 
for some $C_L(p, q, b) < \infty$ if and only if the associated left modulus function 
\begin{equation*}F_L(x) \defeq \norm{s^{1/q}}_{L^q(x, b)} \cdot \norm{m^{-1/p}}_{L^{p'}(0, x)},\ x \in (0,b)
\end{equation*}
satisfies $\sup_{x \in (0, b)} F_L(x) < \infty$. The corresponding embedding 
$V^p_L(0, b; d\sigma) \hookrightarrow L^q(0, b; d\sigma)$ is compact if and only if, additionally, 
\begin{equation*}\lim_{x \to 0^+} F_L(x) = \lim_{x \to b^-} F_L(x) = 0.
\end{equation*}
\end{theorem}
Specialize to $\sigma = \nu_\lambda$ and $dm = dx$ to obtain the following.
\begin{theorem}[Hardy Inequalities for $\lambda > 0$]\label{thm:hardy_lambda_pos}
For $\lambda > 0$ and $d\nu_\lambda(x) = \rho_{\lambda}(x)^{-1} \, dx$:
\begin{enumerate}
  \item If $b < \infty$, then~\eqref{eq:hardy_inequality} (with $d\sigma = d\nu_\lambda$, 
    $d\mu = dx$) holds for all $1 \leq p \leq q < \infty$, with compact embedding 
    $V^p_L((0, b); d\nu_\lambda) \cembed L^q(0, b; d\nu_\lambda)$
    when $1 < p < \infty$.
  \item If $b = \infty$, then~\eqref{eq:hardy_inequality} holds if and only if
    $1 \leq p \leq 2$, $q \leq p'$, and $q < \infty$, with compact embedding 
    $V^{p}_{L}((0,b))\cembed L^q(0, b; d\nu_\lambda)$ if and only if
    $1 < p < 2$ and $2 < q < p'$.
\end{enumerate}
\end{theorem}
\begin{proof}
We verify the conditions of Theorem~\ref{thm:hardy_opic} by computing $F_L(x_0)$ and its limits. 
First, use the change of variable $z = (x_0 + \lambda)/\lambda$ so that $z(x_0) \in (1, (b+\lambda)/\lambda)$ 
and thereby conclude 
\[ \norm{\rho_{\lambda}^{-1/q}}^{q}_{L^q(x_0,b)} = \frac{1}{\lambda}\int_{(x_0+\lambda)/\lambda}^{(b+\lambda)/\lambda}\frac{1}{z\sqrt{z^2-1}}
=\lambda^{-1} \arcsec(z)\Big\vert_{(x_0+\lambda)/\lambda}^{(b+\lambda)/\lambda}\]
We also compute 
\[ \norm{1}_{L^{p'}(0, x_0)} = x_0^{(p-1)/p} = \lambda^{(p-1)/p}(z(x_0) - 1)^{(p-1)/p}\] 
for $p \geq 1$. Altogether, for $1 \leq p < \infty$, 
\[
  F_L(x_0) = \lambda^{(p-1)/p - 1/q} 
  \big[\arcsec\big(\tfrac{b+\lambda}{\lambda}\big) - \arcsec(z(x_0))\big]^{1/q} 
  (z(x_0) - 1)^{(p-1)/p},
\]
For $b < \infty$, $\sup_{x_0 \in (0,b)} F_L(x_0) < \infty,\ \forall\ 1 \leq p < \infty$, 
and only $1 < p < \infty$ gives $\lim_{x_0 \to 0} F_L(x_0) = \lim_{x_0 \to b} 
F_L(x_0) = 0$ (since $\arcsec z(x_0) \to 0$ as $x_0 \to b$ and $z(x_0) \to 1$ 
as $x_0 \to 0$). If $b = \infty$, for $p = 1$ note that $\sup_{x_0 \in (0,\infty)}F_{L}(z(x_0)) < \infty$ as $\arcsec z(x_0) \leq \frac{\pi}{2}$. 
When $p > 1$, we expand $F_L^q(z_0(x_0))$ to leading order as $z_0 \to \infty$ (by expanding $\arcsec$), giving 
\[
  \lim_{z_0 \to \infty} F_L^q(z_0) = (\lambda)^{q\frac{p-1}{p}-1}\lim_{z_0 \to \infty} z_0^{\frac{q(p-1)}{p}-1} =
  \begin{cases} 0 & q(p-1)/p < 1, \\ \lambda^{\frac{q(p-1)}{p}-1} & q(p-1)/p = 1, \\ \infty & q(p-1)/p > 1, \end{cases}.
\]
This information on $F_L$ gives the stated result via Theorem~\ref{thm:hardy_opic}. 
\end{proof}
\begin{theorem}[Hardy Inequalities for $\lambda = 0$]\label{thm:hardy_lambda_zero}
For $\lambda = 0$, we have $d\nu_0(x) = x^{-2} \, dx$. 
\begin{enumerate}
  \item If $b < \infty$, then~\eqref{eq:hardy_inequality} holds if and only if $1 \leq p' \leq q < \infty$, with compact embedding 
    $V^p_0(0, b; d\nu_0) \cembed L^q(0, b; d\nu_0)$ when
    $1 < p < \infty$.
  \item If $b = \infty$ then for any $1 \leq p \leq q < \infty$, \eqref{eq:hardy_inequality} holds if and 
    only if $q = p' = p/(p-1) < \infty$, and the corresponding embedding is not compact. 
\end{enumerate}
\end{theorem}
\begin{proof}
For $p > 1$,
\[
  F_L(x_0) = \big(\tfrac{1}{x_0} - \tfrac{1}{b}\big)^{1/q} \, x_0^{(p-1)/p} 
  = x_0^{\tfrac{p-1}{p} - 1/q} \big(1 - x_0/b\big)^{1/q}
\]
(with $b = \infty$ giving simply $x_0^{\tfrac{p-1}{p} - 1/q}$). For $p = 1$, 
$F_L(x_0) = (x_0^{-1} - b^{-1})^{1/q}$, which is unbounded near $x_0 = 0$. If $b < \infty$ and $p > 1$, $F_L$ is continuous on $(0, b)$ and both endpoint 
limits vanish, giving boundedness and compactness. If $b = \infty$, boundedness of $F_L(x_0) = x_0^{(p-1)/p - 1/q}$ on $(0, \infty)$ 
requires the exponent to vanish, i.e., $q = p/(p-1) = p'$. Inserting this information into
Theorem~\ref{thm:hardy_opic} gives the result. 
\end{proof}
\section{Interpolations}\label{appendix:pointwise_interpolation}
\begin{lemma}[{Parabolic Interpolation, 1D; cf.~\cite[Lemma 7.19]{vazquez}}] \label{lem:holder_interpolation}
For $0 \leq a < b$, $0 < T$, let $\mc{D} = (a,b) \times (0,T)$, suppose
$f \in H^1(\mc{D})$, not identically constant, satisfies $C_x \defeq \norm{\partial_x f}_{L^\infty(0,T;L^2(a,b))} < \infty$.
Then $f$ admits a representative $\tilde{f} \in C^{\frac{1}{4},\frac{1}{2}}(\overline{\mc{D}})$,
and there exists $\mu = \mu(C_x, \norm{\partial_t f}_{L^2(\mc{D})}, b-a) > 0$ such that on any parabolic sub-cylinder
$U = J \times S \subset \mc{D}$ with $\abs{S} \leq \abs{J}^2$ and $\abs{J} \leq \mu$,
\begin{equation}\label{eq:holder_interpolation}
     \abs{\tilde{f}(x_1,t_1) - \tilde{f}(x_2,t_2)} \leq C_x\,\abs{x_1-x_2}^{1/2} + C_t\,\abs{t_2-t_1}^{1/4},
     \quad \forall\,(x_1,t_1),(x_2,t_2) \in U,
\end{equation}
where $C_t = \tfrac{4}{\sqrt{3}}(\norm{\diff{t}f}_{L^2(\mc{D})}C_x)^{1/2}$, and one may take
\[ \mu = \min\Bigl(\tfrac{b-a}{2},\,\tfrac{2(b-a)C_x}{3\norm{\partial_t f}_{L^2(\mc{D})}}\Bigr). \]
\end{lemma}
\begin{proof}
It suffices to prove~\eqref{eq:holder_interpolation} for smooth $f$ on parabolic subcylinders, 
since the estimate for general $f$ follows by approximation, and the conclusion $f \in C^{\frac{1}{4},\frac{1}{2}}(\mc{D})$
is obtained by chaining subcylinders between appropriate points. Take $(t_1,x_1),(t_2,x_2) \in \mc{D}$ and  
introduce $f(t_1,x_2)$ in~\eqref{eq:holder_interpolation} to find 
\begin{align*}
\abs{f(t_1,x_1)-f(t_2,x_2)}
&\leq \abs{f(t_1,x_1)-f(t_1,x_2)} + \abs{f(t_1,x_2)-f(t_2,x_2)} \\
&\defeq I + II.
\end{align*}
By Cauchy-Schwarz, $I \leq \abs{x_2-x_1}^{1/2}C_x$, and it remains to control $II$. 
For $x_2 \in (a,b)$, take $\frac{b-a}{2} \geq r > 0$ to be the one-sided interval $I_r$, either 
$(x_2, x_2+r)$ or $(x_2-r, x_2)$, whichever fits in $(a,b)$. Introduce the average of $f$ over $I_{r}$ 
at each time in $II$ to find the bound, $II \leq (A_1)+(B)+(A_2)$ where $(A_i) := \abs{f(t_i,x_2) - \bint{I_r}f(t_i,\cdot)}$ and
$(B) := \abs{\bint{I_r}[f(t_1,\cdot)-f(t_2,\cdot)]}$. Cauchy-Schwarz on each term gives 
$(A_i) \leq (\bint{I_r}\abs{y-x_2}^{1/2}dy)C_x \leq \frac{2}{3}r^{1/2}C_x$, 
and $(B) \leq r^{-1/2}\abs{t_2-t_1}^{1/2}\norm{\diff{t}}_{L^2(\mc{D})}$. Writing
$E =\abs{t_2-t_1}^{1/2}\norm{\diff{t}f}_{L^2(\mc{D})}$ and $M = \frac{4}{3}C_x$ gives a combined bound
$II \leq E r^{-1/2} + M r^{1/2}$. Optimization in $r$ gives $r_{\ast} = \frac{E}{M}$ and inserting 
$r_{\ast}$ yields $II \leq 2 \sqrt{E M} = \frac{4}{\sqrt{3}}(\norm{\diff{t}f}_{L^2(\mc{D})}C_x)^{1/2}(t_2-t_1)^{1/4}$. 
To enforce that $r_{\ast}$ to be admissible, i.e. $r_{\ast} \leq \frac{b-a}{2}$,
take $(t_1,x_1),(t_2,x_2) \in U = J \times S$ with $S \leq \abs{J}^2$, and, using the formula for $r_{\ast}$, 
impose the requirement $\abs{J} \leq \mu = \min\bigl(\tfrac{b-a}{2},\,\tfrac{2(b-a)C_x}{3\norm{\diff{t}f}_{L^2(\mc{D})}}\bigr)$ to ensure $r_{\ast} \leq \frac{b-a}{2}$ and 
$\abs{J} \leq \frac{b-a}{2}$ and thus the validity of the above estimate for $II$. Combining with the estimate for $I$ gives exactly~\eqref{eq:holder_interpolation} for appropriate subcylinders. 
\end{proof}

\begin{lemma}[Bounded Plus Semiconvex Implies Locally Lipschitz]\label{lem:lip_from_semiconvex}
Let $f \in C^2(B_R(x_0))$ with $|f| \leq M < \infty$ on $B_R(x_0) \subset \bbR^d$ 
satisfy $D^2 f \geq -N$ on $B_R(x_0)$ for some $N > 0$. Then for any compact 
$D \subset B_R(x_0)$ and $\delta < \dist(D, \partial B_R(x_0))$,
\[
  \norm{\grad f}_{L^\infty(D)} \leq \frac{2M}{\delta} + \delta \frac{N}{2}.
\]
\end{lemma}
\begin{proof}
See~\cite[Chapter 2]{CannarsaSinestrari} for general results; we give a direct 
proof here to obtain constants. Fix $x' \in D$ and let $q_{x'}(x) \defeq 
\frac{N}{2} \normsq{x - x'}$, so that $f + q_{x'}$ is convex on $B_R(x_0)$. The 
secant-line characterization of convexity gives, for $\lambda \in (0, 1)$ and 
$x \in B_\delta(x')$,
\[
  \frac{1}{\lambda \norm{x-x'}}\big(f(x' + \lambda(x-x')) - f(x')\big) 
  \leq \frac{f(x) - f(x')}{\norm{x-x'}} + (1 - \lambda) \norm{x-x'} \frac{N}{2}.
\]
Letting $\lambda \to 0$, applying $\abs{f} \leq M$ and $\norm{x - x'} \leq \delta$:
\[
  \frac{x - x'}{\norm{x - x'}} \cdot \grad f(x') \leq \frac{2M}{\delta} + \delta \frac{N}{2}.
\]
Choosing $x$ so that $\frac{x - x'}{\norm{x - x'}} = \frac{\grad f(x')}{\norm{\grad f(x')}}$ and
inserting into the last inequality above gives the conclusion.
\end{proof}
\section{Proof of Claims 1-3 of Proposition~\texorpdfstring{\ref{prop:trunc_gradients_converge}}{3.3}}\label{appendix:dall_aglio}
\newcommand{\yeta}{\bar y_\eta}
We prove Claim 3 first as it differs nontrivially from~\cite[Proposition 6.2]{dallaglio_primary}.
For brevity write $a=1/k$, $y_n=T_ku_n$, $y=T_ku$,
$\yeta=(T_ku)_\eta$, $d_n=y_n-\yeta$, and $E_n=Q_T\cap\{a\leq u_n\leq k\}$.
Thus $a\leq y_n,\,\yeta\leq k$, and $v=\zeta\phi_+(d_n)$ vanishes, together with its spatial derivative, on $\{u_n\leq a\}$.
All limits below are taken first in $n$ and then in $\eta$; $k$ and the nonnegative spatial cutoff $\zeta$ are fixed.

\begin{proof}[Proof of Claim 3]
For choice of $\alpha$ sufficiently large (recall that we take any $\alpha\geq1$ to define $v$ at the start),
\[
 \mc T_3=\fintegral{Q_T}\tilde{u}_n\partial_xu_n\partial_xv
 \geq\fintegral{E_n}\zeta\abssq{\partial_x(d_n)_+}
       +\omega(\eta)+\omega^\eta(n).
\]
Differentiate $v$ and apply $T_k'(s)=\mathbf1_{(a,k)}(s)$ a.e.
to obtain the decomposition
\begin{align*}
 \mc T_3
 &=\fintegral{E_n}\tilde{u}_n
       \abssq{\partial_x(d_n)_+}\phi_+'(d_n)\zeta\\
 &\quad+\fintegral{E_n}\tilde{u}_n
       \partial_x\yeta\,\partial_x(d_n)_+\phi_+'(d_n)\zeta\\
 &\quad-\fintegral{\{u_n>k\}}\tilde{u}_n
       \partial_xu_n\partial_x\yeta\,\phi_+'(k-\yeta)\zeta\\
 &\quad+\fintegral{Q_T}\tilde{u}_n
       \partial_xu_n\phi_+(d_n)\partial_x\zeta
 \defeq A_1+A_2+A_3+A_\zeta.
\end{align*}
There is no lower-tail contribution, since $v$ and $\partial_xv$ vanish
on $\{u_n\leq a\}$.
We find immediately that $A_\zeta$ is vanishing: $\tilde{u}_n\partial_xu_n$ is bounded in $L^2(Q_T)$ by~\eqref{eq:shifted_power_bound},
while $\phi_+(d_n)\partial_x\zeta$ tends to zero strongly in $L^2(Q_T)$ in the successive limits, by dominated convergence. Hence
$A_\zeta=\omega(\eta)+\omega^\eta(n)$.

We next estimate $A_2$ and $A_3$, using Propositions~\ref{prop:uniform_estimates} and~\ref{prop:convergence_from_compactness}.
For $A_2$, write
\[
 A_2=\fintegral{E_n}\tilde{u}_n\partial_x\yeta\,
       \partial_x\phi_+(d_n)\zeta.
\]
Replace $\partial_x\yeta$ by $\partial_xy$, with an error bounded by
$C_k\norm{\partial_x\yeta-\partial_xy}_{L^2(Q_T)}$, since
$0<\tilde{u}_n\leq k+1$ on $E_n$ and the remaining gradient is bounded in $L^2$.
For fixed $\eta$, $\partial_x\phi_+(d_n)$ converges weakly in $L^2$ on the support of $\zeta$ to $\partial_x\phi_+(y-\yeta)$.
Moreover,
$\mathbf1_{E_n}\tilde{u}_n\partial_xy\to y\partial_xy$ strongly in $L^2$,
by dominated convergence: $u_n\to u$ a.e., and $\partial_xy=0$
a.e. on $\{u\notin(a,k)\}$, including the level sets.
Finally, $\partial_x\phi_+(y-\yeta)\to0$ strongly in $L^2$ as $\eta\to\infty$.
Consequently $A_2=\omega(\eta)+\omega^\eta(n)$.

For $A_3$, the flux
\[
 F_{n,k}=\mathbf1_{\{u_n>k\}}\tilde{u}_n\partial_xu_n
 \rightharpoonup \mathbf1_{\{u>k\}}u\partial_xu=:F_k
 \quad\hbox{in }L^2(Q_T)
\]
is uniformly bounded by~\eqref{eq:shifted_power_bound}. Its weak limit is identified by writing $F_{n,k}$ as the spatial derivative of
\[
 \mathbf1_{\{u_n>k\}}
 \left(\frac{u_n^2-k^2}{2}+\frac{u_n-k}{n}\right),
\]
whose local $L^1$ limit follows from strong local $L^2$ convergence of $u_n$.
For fixed $\eta$ we may therefore pass to the limit in $A_3$.
As $\eta\to\infty$, $\phi_+'(k-\yeta)\partial_x\yeta$ converges strongly in $L^2$ to
$\phi_+'(k-y)\partial_xy$, by continuity of the positive-part map in $H^1$ and the smooth exponential on the bounded range.
The latter gradient vanishes on $\{u\geq k\}$, so $A_3=\omega(\eta)+\omega^\eta(n)$.

The remaining contribution is $A_1$. On $E_n$ we have $\tilde{u}_n\geq a=1/k$,
and for $s>0$, $\phi_+'(s)=\alpha e^{\alpha s}\geq\alpha$.
Thus, choosing $\alpha\geq k$ gives
\[
 A_1\geq\frac{\alpha}{k}\fintegral{E_n}
       \abssq{\partial_x(d_n)_+}\zeta
 \geq\fintegral{E_n}\abssq{\partial_x(d_n)_+}\zeta.
\]
Here the gradient of $(d_n)_+$ vanishes a.e. where $d_n\leq0$.
Combining all contributions completes the claim.
\end{proof}

\textbf{Claim 1:}
\[
 \mc T_1=\fintegral{Q_T}\partial_tu_n
       \phi_+(y_n-\yeta)\zeta\,dt\,d\nu_\lambda
       \geq\omega(\eta)+\omega^\eta(n).
\]
\begin{proof}[Proof of Claim 1]
To start, apply the decomposition $s=T_k(s)+G_k(s)$ from~\eqref{eq:identity_cutoff_decomp} to get
\begin{align*}
 \mc T_1
 &=\fintegral{Q_T}\partial_ty_n\phi_+(y_n-\yeta)
                \zeta\,dt\,d\nu_\lambda\\
 &\quad+\fintegral{Q_T}\partial_tG_k(u_n)
       \phi_+(y_n-\yeta)\zeta\,dt\,d\nu_\lambda
       \defeq A_1+A_2.
\end{align*}
We first estimate $A_1$. Introduce the regularized target into the time derivative
and use the ordinary primitive
\[
 \Psi_+(s)=\int_0^s\phi_+(r)\,dr
          =\frac{e^{\alpha s_+}-1}{\alpha}-s_+\geq0.
\]
The chain rule gives
\[
 \partial_ty_n\phi_+(y_n-\yeta)
 =\partial_t\Psi_+(y_n-\yeta)
   +\partial_t\yeta\,\phi_+(y_n-\yeta).
\]
Consequently,
\[
 A_1=\left[\fintegral{\Omega}\Psi_+(y_n-\yeta)\zeta\,d\nu_\lambda\right]_0^T
       +\fintegral{Q_T}\partial_t\yeta\,\phi_+(y_n-\yeta)\zeta\,dt\,d\nu_\lambda.
\]
The terminal term is nonnegative, and the initial term is
$\omega(\eta)+\omega^\eta(n)$, since $T_ku_{0,n}\to T_ku_0$ and
$\yeta(0)=T_ku_{0,j(\eta)}\to T_ku_0$ locally in $L^1(d\nu_\lambda)$,
while $\Psi_+(0)=0$ and $\Psi_+$ is Lipschitz on the bounded range in question.
For fixed $\eta$, dominated convergence and $\partial_t\yeta=\eta(y-\yeta)$ give
\[
 \fintegral{Q_T}\partial_t\yeta\,\phi_+(y_n-\yeta)\zeta\,dt\,d\nu_\lambda
 =\eta\fintegral{Q_T}(y-\yeta)\phi_+(y-\yeta)\zeta\,dt\,d\nu_\lambda
       +\omega^\eta(n)\geq\omega^\eta(n),
\]
because $s\phi_+(s)\geq0$ for every $s\in\bbR$.

To estimate $A_2$, observe that, on the whole real line,
\[
 G_k'(s)=\mathbf1_{(-\infty,a)}(s)+\mathbf1_{(k,\infty)}(s)
 \quad\text{a.e.}
\]
The lower part contributes zero because $\yeta\geq a$ and $\phi_+(a-\yeta)=0$.
Writing $H_n=(u_n-k)_+$, we therefore have
\begin{align*}
 A_2
 &=\fintegral{Q_T}\partial_tH_n\,\phi_+(k-\yeta)
                         \zeta\,dt\,d\nu_\lambda\\
 &=\left[\fintegral{\Omega}H_n\phi_+(k-\yeta)
                         \zeta\,d\nu_\lambda\right]_0^T\\
 &\quad+\fintegral{Q_T}H_n\phi_+'(k-\yeta)\partial_t\yeta
                         \zeta\,dt\,d\nu_\lambda.
\end{align*}
Again the terminal term is nonnegative. The initial term tends to zero in the
successive limits by convergence of the initial data and
$\yeta(0)=T_ku_{0,j(\eta)}$: its limiting integrand is zero since
\[
 (u_0-k)_+\phi_+(k-T_ku_0)=0.
\]
For fixed $\eta$, strong local $L^1$ convergence of $H_n$ shows that the last integral converges to
\[
 \eta\fintegral{Q_T}(u-k)_+\phi_+'(k-\yeta)(y-\yeta)
                         \zeta\,dt\,d\nu_\lambda\geq0,
\]
since $y=k$ wherever $(u-k)_+>0$. The value assigned to $\phi_+'(0)$
does not affect the chain rule, since $\partial_t\yeta=0$ a.e. on $\{\yeta=k\}$.
Thus $A_1$ and $A_2$ are composed of terms which are either nonnegative or asymptotically vanishing, completing Claim 1.
\end{proof}

\textbf{Claim 2:}
\[
 \mc T_2=\fintegral{Q_T}g_n\tilde{u}_n\phi_+(y_n-\yeta)\zeta
       =\omega(n)+\omega^\eta(n)+\omega(\eta).
\]
\begin{proof}[Proof of Claim 2]
Recall that $k$ is fixed. For a compact interval $K\Subset\Omega$ containing the support of $\zeta$, H\"older's inequality gives
\[
 |\mc T_2|
 \leq \linfty{\zeta}\linfty{\rho_\lambda g}
       \norm{\tilde{u}_n}_{L^2((0,T)\times K;dt\times d\nu_\lambda)}
       \norm{\phi_+(y_n-\yeta)}_{L^2((0,T)\times K;dt\times d\nu_\lambda)}.
\]
The first norm is uniformly bounded by Proposition~\ref{prop:uniform_estimates} and $\nu_\lambda(K)<\infty$.
The second tends to zero, first in $n$ and then in $\eta$, by dominated convergence on the bounded range $a\leq y_n,\,\yeta\leq k$. This completes the claim.
\end{proof}

\section*{Acknowledgments}
We thank Professor Juan Luis V\'azquez for enlightening discussions
on limits of porous medium equations during his visit to Austin,
Texas. We thank Professor Stefania Patrizi for her insightful
instruction on regularity theory of parabolic equations.
The authors acknowledge the support and funding of the UT Austin
Analysis RTG, NSF grants DMS-2009736, DMS-RTG-1840314, and
DMS-2408263, and Department of Energy grant DOE-DESC0016283,
Project Simulation Center for Runaway Electron Avoidance and
Mitigation.

\section*{Statements and Declarations}

\noindent\textbf{Competing interests.}
The authors declare that they have no competing interests.

\medskip\noindent\textbf{Generative AI use.}
The authors used Claude and ChatGPT to assist with writing and
editing prose, LaTeX preparation, proof checking and editing,
and idea exploration. All results and proofs were developed
originally and directly by the authors, who assume full
responsibility for the manuscript.

\medskip\noindent\textbf{Data availability.}
No datasets were generated or analyzed during this study.

\bibliographystyle{plain} 
\bibliography{references}

\end{document}